\RequirePackage{fix-cm}
\RequirePackage{amsmath}

\documentclass[]{article}
\usepackage{latexsym}
\usepackage{mathtools}
\mathtoolsset{showonlyrefs}
\usepackage{listings}
\usepackage{xcolor}
\usepackage{fullpage}
\usepackage{braket}
\usepackage{graphicx,epstopdf,amssymb,amsfonts,mathrsfs,enumitem,color,hyperref,url,bm}
\usepackage{caption}
\usepackage{subcaption}

\usepackage{algorithm,algorithmic}

\newcommand{\wt}{\widetilde}

\newcommand{\f}{\mathbb}
\newcommand{\cu}{\subseteq}

\newcommand{\mc}{\mathcal}

\newcommand{\IMF}{\textrm{IMF}}

\DeclareMathOperator{\diag}{diag}

\newtheorem{lemma}{Lemma}
\newtheorem{theorem}{Theorem}
\newtheorem{proposition}{Proposition}
\newtheorem{corollary}{Corollary}
\newtheorem{definition}{Definition}

\begin{document}

\title{Stabilization and Variations to the Adaptive Local Iterative Filtering Algorithm: the Fast Resampled Iterative Filtering Method}

\author{Giovanni Barbarino\thanks{Faculté Polytechnique de Mons, Université de Mons, Mons, Belgium. giovanni.barbarino@gmail.com}, Antonio Cicone\thanks{DISIM, Universit\`a degli Studi dell'Aquila, L'Aquila, Italy, and Istituto di Astrofisica e Planetologia Spaziali, INAF, Roma, Italy, and
Istituto Nazionale di Geofisica e Vulcanologia, Roma, Italy. antonio.cicone@univaq.it}}

\maketitle

\begin{abstract}
Non-stationary signals are ubiquitous in real life. Many techniques have been proposed in the last decades which allow decomposing multi-component signals into simple oscillatory mono-components, like the groundbreaking Empirical Mode Decomposition technique and the Iterative Filtering method. When a signal contains mono-components that have rapid varying instantaneous frequencies, we can think, for instance, to chirps or whistles, it becomes particularly hard for most techniques to properly factor out these components.
The Adaptive Local Iterative Filtering technique has recently gained interest in many applied fields of research for being able to deal with non-stationary signals presenting amplitude and frequency modulation.
In this work, we address the open question of how to guarantee a priori convergence of this technique, and propose two new algorithms.
The first method, called Stable Adaptive Local Iterative Filtering, is a stabilized version of the Adaptive Local Iterative Filtering that we prove to be always convergent. The stability, however, comes at the cost of a higher complexity in the calculations. The second technique, called Resampled Iterative Filtering, is a new generalization of the Iterative Filtering method. We prove that Resampled Iterative Filtering is guaranteed to converge a priori for any kind of signal. Furthermore, in the discrete setting, by leveraging on the mathematical properties of the matrices involved, we show that its calculations can be accelerated drastically. Finally, we present some artificial and real-life examples to show the powerfulness and performance of the proposed methods.
\end{abstract}

\section{Introduction}\label{sec:introduction}

The analysis and decomposition of non-stationary signals is an active research direction in both Mathematics and Signal Processing. In the last decades many new techniques have been proposed. Among them, the Iterative Filtering (IF) algorithm \cite{lin2009iterative} was proposed a decade ago as an alternative technique to the celebrated Empirical Mode Decomposition (EMD) technique \cite{huang1998empirical} and its variants \cite{wu2009ensemble,yeh2010complementary,torres2011complete,zheng2014partly}. The EMD and its variants, in fact, were missing a rigorous mathematical analysis, due to the usage of a number of heuristic and ad hoc elements. Some results have been presented in the literature \cite{huang2009convergence,ur2011filter,huang2014introduction}, but a complete rigorous mathematical analysis is still missing nowadays.

The EMD-like methods are based on the iterative computation of the signal moving average via envelopes connecting its extrema. The computation of the signal moving average allows to split the signal itself into a small number of simple oscillatory components, called Intrinsic Mode Functions (IMFs), which are separated in frequencies and almost uncorrelated \cite{cicone2021IMFogram}. The IF method has been developed following the same structure of EMD, but with a key difference: the moving average is now obtained through an iterated convolutional filtering operation on the signal, with the aim to single out all its non-stationary components, starting from the highest frequency one.

The IF algorithm structure allowed to develop a complete mathematical analysis of this method \cite{cicone2020iterative,cicone2016adaptive,cicone2020numerical}. On the other side, this method is ``rigid'' in the sense that it allows to extract only IMFs which are amplitudes modulated, but almost stationary in frequencies. This is a clear limitation if the signal contains chirps or whistles, that are components with quickly changing instantaneous frequencies. For this reason in \cite{cicone2016adaptive} the authors proposed a generalization of IF called Adaptive Local IF (ALIF). ALIF does not suffer anymore of the rigidity of IF in extracting IMFs containing rapidly varying instantaneous frequencies. However, this new technique loses most of mathematical background of IF. Even if the algorithm gained visibility since its introduction five years ago, we can mention here, for instance, \cite{An2017,%
	an2016application,%
	an2016wind,%
	an2017vibration,%
	an2016demodulation,%
	kim2016multiscale,%
	li2018entropy,%
	mitiche2018classification,%
	cicone2017Geophysics,sharma2017automatic,yang2017oscillation}, an initial mathematical analysis has been only recently developed \cite{cicone2019spectral,cicone2020convergence}, and much more study on extensions, variations and stabilization methods is currently ongoing, see, for instance, \cite{barbarino2021conjectures}.

Due to the missing theoretical background of the ALIF method, in this paper we introduce two new algorithms for which such analysis is possible. The first, called Stable ALIF (SALIF) method, is always convergent, even in presence of noise, but it presents an increased computational cost with respect to ALIF. The second, called Resampled IF method (RIF), is actually a modification of the IF algorithm, that preserves IF convergence property, but, at the same time, presents the same flexibility as ALIF. Furthermore RIF method can be made, in the discrete case, highly computational efficient via the FFT computation of the convolutions, in what is called the Fast Resampled IF method (FRIF).

The rest of this paper is structured as follows. Section \ref{sec:IF_and_ALIF} is a review of the IF and ALIF methods, and introduces the new SALIF method. Here we compare their features, stressing their strength and weaknesses. Section \ref{sec:RIF} is dedicated to the RIF algorithm, its analysis, properties and acceleration via FFT, in what is called FRIF technique. In this section we show how RIF combines the convergence and stability of IF with the flexibility of ALIF, and how it can be made computationally efficient. In Section \ref{sec:NumericalExamples} we compare those algorithms on artificial and real data, reporting the efficiency and accuracy of each method. Eventually, in Section \ref{sec:Outlook}, we draw conclusions and  suggest future lines of research.

\section{Iterative Filtering based Methods}\label{sec:IF_and_ALIF}

Throughout this document, a signal is intended to be a real function $g:\f R\to \f R$, and we study its behaviour in the reference interval $[0,1]$. Outside this interval the signal is usually not known, and so we have to impose some boundary conditions, discussed for example in \cite{cicone2020study} and \cite{stallone2020new}. In particular, in \cite{stallone2020new}, the authors show how any signal can be pre-extended and made periodical at the boundaries. Therefore, from now on, for simplicity and without losing generality, we will assume that the signals to be decomposed are always periodical at the boundaries.

The Iterative Filtering (IF) methods mimics the EMD algorithm in the application of a moving average that captures the main trend of the signal, and allows us to decompose it into simple IMF components. If we call $\mc L(g)$ the moving average, then both EMD and IF algorithms extract the first IMF as
\begin{equation}\label{eq:first_IMF}
	\mc S(g) = g - \mc L(g), \qquad IMF_1 = \lim_{m\to\infty} \mc S^m(g).
\end{equation}
Repeating iteratively the same procedure on $r = g - IMF_1$, we can extract all the IMFs until $r$ becomes a trend signal, meaning that it possesses at most two extrema.

The difference between these two algorithms is that, while for EMD the moving average operator $\mc L(g)$ is changing at each iteration and depends completely on the shape of a given signal, in IF $\mc L(g)$ can be rewritten as the convolution of $g$ with what is called a filter $k$.
Here a filter $k$ is an even, nonnegative, bounded and measurable real function  with compact support and unit mass, meaning $\int_{\f R} k(z)\,\,{{\rm d}}z=1$.

A generalization of the IF method is called Adaptive Local Iterative Filtering (ALIF), and utilizes the convolution with a family of filters $k_x(z)$ as moving average, whose support is in $[-\ell(x),\ell(x)]$, i.e. it varies with $x$. Therefore, the moving average computation operator can be written as
\begin{equation}\label{eq:moving_average_kx}
	\mc L(g)(x) = \int_{0}^{1} g(z)k_x(x-z) \,\,{{\rm d}}z,
\end{equation}

Following \cite{cicone2016adaptive}, we can rewrite the same
expression as
\begin{equation}\label{eq:moving_average_txz}
\mc L(g)(x) = \int_{-L}^L g(x + t(x,z) )k(z) \,\,{{\rm d}}z,
\end{equation}
where $k(z)$ is a filter with constant support in $[-L,L]$ and $t(x,z)$ is a measurable function. In the following subsections we report the most common choice for the filters, and a description of the resulting method.

\subsection{Linear ALIF}

When we talk about the ALIF method, we usually refer to  Linear ALIF. After having fixed a filter $k(z)$ with support $[-1,1]$, and a  positive ``length'' function $\ell(x)$, then the linear ALIF method prescribes
\[
k_x(z) := k\left(  \frac z{\ell(x)}\right)\frac 1{\ell(x)}\text{ in \eqref{eq:moving_average_kx}}, \quad \text{ or equivalently } \quad
t(x,z) = -\ell(x)z\text{ in \eqref{eq:moving_average_txz}}.
\]
Notice that $k_x(z)$ is a filter with support $[-\ell(x), \ell(x)]$ for every $x\in [0,1]$.

Given now a signal $g(z)$, we can compute a length function $\ell(x)$, that usually depends on the relative positions of local extrema in $g(z)$ if the signal does not contain noise, and apply the iteration in \eqref{eq:first_IMF} with the appropriate filter.
\begin{equation}\label{eq:first_IMF_ALIF}
\mc S(g)(x) = g(x) - \int_0^1 g(z)k\left(  \frac {x-z}{\ell(x)}\right)\frac {{{\rm d}}z}{\ell(x)} , \qquad IMF_1 = \lim_{m\to\infty} \mc S^m(g).
\end{equation}
Repeating iteratively the same procedure on $r = g - IMF_1$, we obtain a decomposition of the signal $g$ into IMFs. Notice that $\ell(x)$ changes after we identify each different IMFs.
Here we report the resulting algorithm.

\begin{algorithm}
	\caption{\textbf{(ALIF Algorithm)} ${\rm IMFs=ALIF}(g)$}
	\begin{algorithmic}
		\STATE IMFs = $\left\{\right\}$
		\STATE initialize the remaining signal $r=g$
		\WHILE{the number of extrema of $r$ is $\geq 2$}
		\STATE for each $x\in[0,1]$ compute the length function $\ell(x)$, depending on $r$
		\STATE $g_1=r$
		\STATE $m=1$
		\WHILE{the stopping criterion is not satisfied}
		\STATE $g_{m+1}=g_m-\int_0^1 g_m(z)k\left(  \frac {x-z}{\ell(x)}\right)\frac {{{\rm d}}z}{\ell(x)}$
		\STATE $m=m+1$
		\ENDWHILE
		\STATE ${\rm IMFs}={\rm IMFs}\cup\{g_m\}$
		\STATE $r=r-g_m$
		\ENDWHILE
		%\STATE ${\rm IMFs=IMFs}\cup\{r\}$
	\end{algorithmic}
	\label{alg:ALIF}
\end{algorithm}

The operation $\mc S(g) = g - \mc L(g)$ is designed to catch the fluctuation part of the signal, that usually presents high frequency. The operation is iterated until a stopping criterion is satisfied, usually regarding the norm of the difference $g_{m+1}-g_m$, or the number of iterations themselves. For more details on the stopping criterion, we refer the interested reader to \cite{cicone2016adaptive,lin2009iterative}.
The IMFs are thus extracted from the signal until it becomes just a trend signal with $2$ or less extrema. Since the sum of all the IMFs and the trend signal returns the original signal, it can effectively be called a decomposition.

Regarding the length function $\ell(x)$ identification in signals containing noise, we observe that it is always possible to first run a time-frequency representation (TFR) algorithm, see \cite{WU2020current} for a comprehensive review of modern TFR techniques, and then use the acquired information to designed the optimal $\ell(x)$. This procedure is really important for ALIF algorithm, but it is also a research topic per se. This is why, from now on, we  assume that the length function can be computed accurately and we postpone the analysis of how actually compute it to a future work.

Conceptually, the ALIF method separates non-stationary components of the signal, even with varying amplitudes, starting from the highest frequencies.
For example, on real data, the method first extract high frequency noise IMFs, and then starts to produce clean components. The main feature of the produced IMFs is that their instantaneous frequencies are pointwise sorted in decreasing order. In formulae, if $F_k(x)$ is the instantaneous frequency of the $k$-th IMF at the point $x\in (0,1)$, we have that
\[
F_1(x) > F_2(x) > F_3(x) > \dots \qquad \forall x\in (0,1).
\]
The method, albeit being very powerful and having already been utilized in a variety of applications, still lacks a theoretical analysis proving the convergence of \eqref{eq:first_IMF_ALIF}, except in a few notorious cases \cite{barbarino2021conjectures,cicone2019spectral,cicone2016adaptive}.
In the next sections, we report some of the available convergence results for the discrete version of the algorithm.

\subsection{Discrete ALIF and Stabilization}

Usually in a discrete setting, a signal is given as a vector of sampled values $\bm g = [g_0\,\, g_1\,\,\dots\,\,g_{n-1}]$ where $g_i = g(x_i)$ and $x_i = i/n$ for $i\in \f Z$. As a consequence, one can discretize the relation \eqref{eq:moving_average_kx} with a simple quadrature formula.
\begin{equation}\label{eq:moving_average_dkx}
	\mc L(g)(x_i) = \int_{0}^{1} g(z)k_{x_i}(x_i-z) \,\,{{\rm d}}z \approx \frac 1{n}\sum_{j=0}^{n-1} g_jk_{x_i}(x_i-x_j)
\end{equation}
In turn, this lets us write the sampling vector of $\mc S(g)$, called $\bm{\wt g}$, as a matrix-vector multiplication. In fact, if we assume all the indexes start from zero,
\begin{equation}\label{eq:matrix_K}
\mc S(g)(x_i) = g_i - \mc L(g)(x_i)
\implies
\bm{\wt g} = (I - K)\bm g,\qquad
K_{i,j} = \hspace{-1pt}\left[\frac 1{n}\,k_{x_i}(x_i-x_j)\right]_{i,j=0}^{n-1}\hspace{-2pt}.
\end{equation}
In the Linear ALIF paradigm, we fix a filter $k(z)$ and choose a length function $\ell(x)$ depending on the signal, to produce our family of filters $k_{x}(z) = k(z/\ell(x))/\ell(x)$. The resulting algorithm is reported here.
\begin{algorithm}\label{alg:dALIF}
	\caption{\textbf{(Discrete ALIF Algorithm)} $\bm{IMFs}={\rm ALIF}(\bm g)$}
	\begin{algorithmic}
		\STATE $\bm{IMFs} = \left\{\right\}$
		\STATE initialize the remaining signal $\bm r=\bm g$
		\WHILE{the number of extrema of $\bm r$ is $\geq 2$}
		\STATE  compute $\ell(x)$ and the matrix $K$
		\STATE $\bm g_1=\bm r$
		\STATE $m=1$
		\WHILE{the stopping criterion is not satisfied}
		\STATE $\bm g_{m+1}=(I-K)\bm g_m$
		\STATE $m=m+1$
		\ENDWHILE
		\STATE $\bm{IMFs}=\bm{IMFs}\cup\{\bm g_m\}$
		\STATE $\bm r=\bm r-\bm g_m$
		\ENDWHILE
		%\STATE $\bm{IMFs}=\bm{IMFs}\cup\{\bm r\}$
	\end{algorithmic}
	\label{alg:ALIFd}
\end{algorithm}

From the algorithm, it is evident that the convergence of the internal loop only depends on the spectral properties of the matrix $I-K$. In fact, since $\bm g_{m+1} = (I-K)^m\bm g_1$, we find that a necessary condition for the convergence is
\begin{equation}\label{eq:convergence_spectral_condition_LALIF}
	|1-\lambda_i(K)| \le 1, \qquad \forall i=1,\dots,n.
\end{equation}
If the zero eigenvalue of $K$ has equal geometric and algebraic multiplicities, and it is the only eigenvalue for which $|1-\lambda_i(K)| = 1$,
 then the condition is also sufficient. From the same analysis, one can notice that the algorithm actually produces a projection of the signal $\bm g$ on the approximated null space of $K$.\\

Notice that if $A$ is an invertible matrix, then $Null(AK)=Null(K)$, meaning that substituting $AK$ in the algorithm  doesn't significantly change the output. As a consequence, we can always suppose that $K$ is a stochastic matrix, so that we don't have to worry about large $\lambda_i(K)$. Nonetheless, we still cannot assert the absence of negative or complex eigenvalues which are not fulfilling the relation \eqref{eq:convergence_spectral_condition_LALIF}. Recent studies \cite{cicone2019spectral,barbarino2021conjectures} show that for large $n$ and continuous functions $k(z)$, $\ell(x)$, almost all eigenvalues of the matrix $K$ are real and nonnegative, but it is still not enough to establish the convergence of the method. Moreover, it has been ascertained experimentally that such cases may arise, especially whit a fast changing function $\ell(x)$.

A simple way to stabilize the method is to choose $A=\|K\|^{-2}K^T$, so that $AK$ is a nonnegative matrix that is also positive semidefinite, with all the eigenvalues bounded by $1$. Notice that we can also use $c^{-2}$ instead of $\|K\|^{-2}$, where  $c = \max_j \sum_i K_{i,j}$, or in general any constant satisfying $\|cK\|\le 1$. We call the resulting method Stable ALIF (SALIF).
\begin{algorithm}
	\caption{\textbf{(Stable Discrete ALIF Algorithm)} $\bm{IMFs}={\rm SALIF}(\bm g)$}
	\begin{algorithmic}
		\STATE $\bm{IMFs} = \left\{\right\}$
		\STATE initialize the remaining signal $\bm r=\bm g$
		\WHILE{the number of extrema of $\bm r$ is $\geq 2$}
		\STATE  compute $\ell(x)$ and the matrix $K$
		\STATE $\bm g_1=\bm r$
		\STATE $m=1$
		\WHILE{the stopping criterion is not satisfied}
		\STATE $\bm g_{m+1}=(I - K^TK)\bm g_m$
		\STATE $m=m+1$
		\ENDWHILE
		\STATE $\bm{IMFs}=\bm{IMFs}\cup\{\bm g_m\}$
		\STATE $\bm r=\bm r-\bm g_m$
		\ENDWHILE
		%\STATE $\bm{IMFs}=\bm{IMFs}\cup\{\bm r\}$
	\end{algorithmic}
	\label{alg:SALIF}
\end{algorithm}

The method is called stable since a perturbation of the matrix $K$ does not prevent the convergence of the inner loop. Moreover, we will show in the experiments, ref. Section \ref{sec:NumericalExamples}, that SALIF is able to produce more accurate solutions than the other methods.

The algorithm, though, comes with an increased computational cost with respect to ALIF, mainly due to two factors.
\begin{itemize}
	\item The iterative step in the SALIF algorithm $\bm g_{m+1}=(I - K^TK)\bm g_m$ takes at least double the time with respect to the respective step in the ALIF algorithm.
	Since the number of iterations is usually much smaller than $n$, even computing $K^TK$ beforehand does not improve the speed.
	\item The order of the smallest eigenvalues of $K^TK$ is approximately the square of the smallest ones in $K$.  The algorithm thus requires more iterations to attain the same accuracy of ALIF, since it must separate eigenspaces that are now closer.
\end{itemize}
A different way to stabilize the method is to take $\ell(x)$ constant, producing a much faster algorithm, i.e. the IF algorithm, whose spectrum of application is though more limited.

\subsection{IF and Discrete IF}

When we talk about the IF method, we refer to the linear ALIF method with constant length function $\ell(x) = L$ in \eqref{eq:first_IMF_ALIF}, or equivalently, where $t(x,z) = -Lz$ in \eqref{eq:moving_average_txz}.
The IF method only separates IMF components of the signal which are amplitude modulated, but quasi-stationary in frequency, starting from the highest frequencies. Nevertheless, it has been proved \cite{huang2009convergence} that in this case the iterations \eqref{eq:first_IMF_ALIF} always converge whenever  $k(z)$ is a filter with nonnegative Fourier transform. The condition is satisfied, for example, by $k = \omega \star \omega$, where $\omega$ is a generic filter and $\star$ is the convolution operator.

In the discrete setting, the IF algorithm has the advantage of a fast implementation based on FFT, in what is called Fast Iterative Filtering (FIF), and an advanced theoretical analysis \cite{cicone2020iterative,cicone2020numerical}. Recall that we only know the signal $g(x)$ on the interval $[0,1]$, so we can always suppose that the original signal is $1$-periodic (for example, by reflecting the signal on both sides and making it decay \cite{stallone2020new}). We can thus rewrite the moving average \eqref{eq:moving_average_kx} as
\begin{equation}\label{eq:moving_average_IF}
\mc L(g)(x) = \int_{\f R} g(z)k\left(\frac{x-z}{L}\right)\frac {{{\rm d}}z}L ,
\qquad
\mc S(g) = g - \mc L(g)
\end{equation}
and discretize it on a regular grid $x_i = i/n$ of $[0,1]$. Here, the integral is always well-defined, since the filter has compact support. Moreover $L$ is inversely proportional to the target frequency of the extracted IMF, and $L\ge 1/2$ usually indicates that we already have a trend signal $g$, so we always suppose $1/L> 2$.

Following the same steps as in the ALIF algorithm, we find
\[
\mc L(g)(x_i)
\approx \frac 1{nL}\sum_{j\in \f Z}g_j
k\left(\frac{x_{i-j}}{L}\right)
\]
Notice that the above formula can be expressed through a Hermitian circulant matrix $K$ with first row
\[
\frac{1}{nM}
\left[
k(0), \, k\left( \frac{1}{nL}  \right)
, \dots
, \, k\left( \frac{s}{nL}  \right)
, \, 0
\dots
\, 0
, \, k\left( \frac{s}{nL}  \right)
, \dots
, \, k\left( \frac{1}{nL}  \right)
\right]
\]
where $s = \lfloor nL \rfloor < \lfloor n/2 \rfloor$.
The sampling vector $\bm{\wt g}$ of $\mc S(g)$ on the points $x_0,x_1,\dots, x_{n-1}$, can be thus rewritten as a matrix-vector multiplication
\begin{equation}\label{eq:matrix_K_IF}
\mc S(g)(x_i) = g_i - \mc L(g)(x_i)
\implies
\bm{\wt g} = (I - K)\bm g.
\end{equation}
The resulting algorithm is thus the same as Algorithm 2, but where $K$ is Hermitian and circulant.
The IF method is consequently much faster than the ALIF algorithm since the multiplication $(I-K)\bm g_m$ can be performed very efficiently through an FFT. Actually, in \cite{cicone2020numerical} we can find an even faster implementation, the so called FIF algorithm, and the proof that $K$ is also positive semidefinite.

Keeping in mind that, as in ALIF, we can always multiply $K$ by a diagonal matrix and make it stochastic, we have the following result.

\begin{lemma}[{\cite[Theorem 1, Corollary 3]{cicone2020numerical}}]\label{res:IF_convergence}
Given a double-convoluted filter $k = \omega\star \omega$, then for the IF operator $\mc S(\cdot)$ in \eqref{eq:moving_average_IF}  the limit
\[
\lim_{m\to\infty} \mc S^m(g)
\]
converges for any function $g(x)$.
Moreover, if $L\le 1/2$, then for the IF matrix $K$ in \eqref{eq:matrix_K_IF} the limit
\[
\lim_{m\to\infty}
(I-K)^m\bm g
\]
converges for any vector $\bm g$.
\end{lemma}

\noindent To summarize, we have
\begin{itemize}
	\item the IF and FIF algorithms always converge and are very fast, but cannot capture non-stationary components with quickly varying frequencies,
	\item the ALIF algorithm is enough flexible to extract fully non-stationary components, but its convergence is not guaranteed,
	\item the SALIF algorithm is always convergent and it has an output which is more accurate than the ALIF algorithm one, but it is very slow.
\end{itemize}
In the next section, we show how to design an alternative method, that is flexible enough to perform non-stationary analysis on the signals, but at the same time fast and provably convergent.

\section{Resampled Iterative Filtering}\label{sec:RIF}

The linear ALIF method makes use of a length function $\ell(x)$ to locally stretch a fixed filter $k(z)$ so that the convolution with the signal $g(z)$ smooths out the high oscillatory behaviour. The idea behind the Resampled Iterative Filtering (RIF) algorithm is  to set a fixed length for $k(z)$ and instead modify the signal through a global resampling function. In a sense, we want to locally stretch the signal, making the component of higher frequency approximately stationary, so that we are able to identify it through the fast IF algorithm.

Given a resampling function
$G: [0,M]\to [0,1]$ that is increasing and regular enough,  the moving average for the RIF method will coincide with the IF one applied on $g\circ G$, as
\begin{equation}\label{eq:moving_average_RIF}
\mc L(g\circ G)(y) = \int_{\f R} g(G(z))k(y-z)\,\, {{\rm d}}z,
\end{equation}
where we assume that the resampled signal is $M$-periodic. If we consider the first-order expansion of $G(x)$, and after a change of variable $x=G(y)$, $r = zG'(y)$, we have
\begin{align}\label{eq:LALIF_as_first_order_RIF}
\nonumber \int_{\f R} g(G(z))\,k(y-z)\,\, {{\rm d}}z &
= \int_{\f R} g(G(y-z))\, k(z)\,\, {\rm d}z\\
\nonumber
&
\approx  \int_{\f R} g(G(y)-zG'(y))\, k(z)\,\, {\rm d}z
\\ \nonumber &
= \int_{\f R} g(x-r)\, k\left( \frac{r}{G'(G^{-1}(x))}\right)
\frac{{\rm d}r}{G'(G^{-1}(x))}\\
&
=  \int_{\f R} g(r) \,k\left( \frac{x-r}{G'(G^{-1}(x))}\right)
\frac{{\rm d}r}{G'(G^{-1}(x))}
\end{align}
that is analogous to the linear ALIF moving average in \eqref{eq:first_IMF_ALIF}, where, equivalently,
\begin{equation}\label{eq:resampling_from_length}
	\ell(x) = G'(G^{-1}(x)) \quad \text{ or }\quad
	G^{-1} (x) = \int_0^x \frac 1{\ell(t)} \,{{\rm d}}t.
\end{equation}
%Actually, with a similar analysis, we can find that any method prescribing a degree $p$ polynomial in $z$ as $t(x,z)$ in \eqref{eq:moving_average_txz} can be seen as a $p$-th order approximation of RIF.
With \eqref{eq:resampling_from_length}, now we have a way to derive the resampling function from the length $\ell(x)$. The full RIF algorithm is thus reported as Algorithm \ref{alg:RIFc}.

\begin{algorithm}
	\caption{\textbf{(Resampled IF Algorithm)} ${\rm IMFs=RIF}(g)$}
	\begin{algorithmic}
		\STATE IMFs = $\left\{\right\}$
		\STATE initialize the remaining signal $r=g$
		\WHILE{the number of extrema of $r$ is $\geq 2$}
		\STATE compute $\ell(x)$  and derive the resampling $G(y), G^{-1}(x)$ and the resampled signal $h = r\circ G$
		\STATE $h_1=h$
		\STATE $m=1$
		\WHILE{the stopping criterion is not satisfied}
		\STATE $h_{m+1}=h_m-\int_{\mathbb R}h_m(y)k(x-y){\rm d}y$
		\STATE $m=m+1$
		\ENDWHILE
		\STATE ${\rm IMFs}={\rm IMFs}\cup\{h_m\circ G^{-1}\}$
		\STATE $r=r-h_m\circ G^{-1}$
		\ENDWHILE
		%\STATE ${\rm IMFs=IMFs}\cup\{r\}$
	\end{algorithmic}
	\label{alg:RIFc}
\end{algorithm}

From the algorithm it is evident that, after the resampling, the steps are the same of the IF algorithm. In fact, we always extract almost stationary IMFs from the resampled signal, and then we operate the inverse sampling to obtain the respective IMFs for the original signal. Moreover, we point out that $G(x)$ depends on $\ell(x)$, so it must be computed every time we want to extract a new component.

This observation is also enough to show that the internal loop always converge to some IMF. In the next section we see how these properties carry to the discrete case.

Notice that RIF is actually a particular ALIF method, since
\[
	\int_{\f R} g(G(x-z))\, k(z)\,\, {\rm d}z
	=
	\int_{\f R} g(x + [G(x-z)-x])\, k(z)\,\, {\rm d}z
\]
with $t(x,y) = G(x-z)-x$ in \eqref{eq:moving_average_txz}.
	
	From the relations \eqref{eq:LALIF_as_first_order_RIF}, we can say that Linear ALIF is a first-order approximation of RIF, and since RIF is a convergent method, we could ask whether it produces the same output as Linear ALIF. The answer is provided in the following
\begin{theorem}
The RIF method produces the same output as the Linear ALIF algorithm only when $\ell(x)$ is a constant function, i.e. when the linear ALIF algorithm reduces to IF.
\end{theorem}
The proof follows from the observation that the derivation of equation \ref{eq:LALIF_as_first_order_RIF} holds true only if
\[
	G(y-z) = G(y) - G'(y)z, \qquad \forall y,z
\]
meaning that $G'(z) = \ell(G(z))$ is constant for every $z$.

\subsection{Fast Resampled Iterative Filtering}

First of all we review how to possibly implement a discrete version of RIF. One way is by discretizing the IF moving average on the resampled signal, as in
$$ g_{m+1}(G(x)) = g_m(G(x)) - \int_{\f R} g_m(G(y))\, k(x-y) {\rm d}y.$$
Notice that $h_m = g_m\circ G$ has domain $[0,M]$, so we need to discretize it on the regular grid $x_i :=Mi/n$ for $i=0,\dots,n-1$. Recall that $M = \int_0^1 \ell(x)^{-1} \,{{\rm d}}x$ and that in the IF algorithm, a constant $L\ge 1/2$ indicates that $g(x)$ is already a trend signal. This shows that we can safely assume $\ell(x) < 1/2$ and thus $M>2$.

We now extend the signal cyclically on the real line, meaning that $h_m(sM + x):= h_m(x)$ for every $s\in\f Z$ and every $x\in [0,M)$. The quadrature rule on the discretization points yields
\[ h_{m+1}(x_i) \approx h_m(x_i) - \frac{M}{n}\sum_{j\in\mathbb Z}h_m(x_j) k(x_i - x_j). \]
Notice that the above formula coincides with the IF moving average with length $L = 1/M$, and can be expressed through a Hermitian circulant matrix $K$ with first row
\[
\bm k_1=\frac{M}{n}
\left[
k(0), \, k\left( \frac{M}{n}  \right)
, \dots
, \, k\left( s\frac{M}{n}  \right)
, \, 0
\dots
\, 0
, \, k\left( s\frac{M}{n}  \right)
, \dots
, \, k\left( \frac{M}{n}  \right)
\right]
\]
where $s = \lfloor n/M\rfloor$.
The moving average thus becomes
$$
\bm{h}_{m+1} = (I- K)\bm h_m
$$
where $I-K$ is still a Hermitian and circulant matrix, so that the matrix vector multiplication can be performed efficiently through a FFT. In particular,
\[
\bm h_{m+1} = \textrm{iDFT}\left((1 - \textrm{DFT}(\bm k_1))  \circ \textrm{DFT}(\bm h_m)\right),
\]
where $\circ$ stands for the Hadamard (or element-wise) product between vectors, and DFT, iDFT stand for Discrete Fourier Transform and its inverse, respectively. Moreover, since
\[
\textrm{DFT}(\bm h_{m+1}) =  (1 - \textrm{DFT}(\bm k_1))  \circ \textrm{DFT}(\bm h_m)
\]
and since the stopping criterion can be checked on $\textrm{DFT}(\bm h_m)$, we can further accelerate the method by computing the DFTs on $\bm h_1$ and $\bm k_1$ and the iDFT outside the loop, thus avoiding iterated computations of Fourier transforms.

 The resulting method is reported in Algorithm \ref{alg:FRIF}.

\begin{algorithm}
	\caption{\textbf{(Fast Resampled Iterative Filtering)} $\bm{IMFs}={\rm FRIF}(\bm g)$}
	\begin{algorithmic}
		\STATE IMFs = $\left\{\right\}$
		\STATE initialize the remaining signal $\bm r=\bm g$
		\WHILE{the number of extrema of $\bm r$ is $\geq 2$}
		\STATE compute $\ell(x)$, the resampling functions
		$G^{-1}(x) = \int \ell(t)^{-1}\,{{\rm d}}t$,
		the constant $M = G^{-1}(1)$ and the matrix $K$
		\STATE compute the vector $\bm h$ through interpolation of $\bm r$ on the points $G(y_i)$ where $y_i = Mi/n$
		\STATE $\bm h_1=\bm h$
		\STATE $\bm{\hat h_1} = \textrm{DFT}(\bm h_1)$,
		$\bm{\hat k_1} = 1-\textrm{DFT}(\bm k_1)$
		\STATE $m=1$
		\WHILE{the stopping criterion is not satisfied}
		\STATE
		$\bm {\hat h_{m+1}}= \bm{\hat k_1} \circ \bm{\hat h_m}$
		\STATE $m=m+1$
		\ENDWHILE
		\STATE $\bm h_m = \textrm{iDFT}(\bm {\hat h_m})$
		\STATE compute the vector $\bm I$ through interpolation of $\bm h_m$ on the points $G^{-1}(y_i)$
		where $y_i = i/n$
		\STATE ${\rm IMFs}={\rm IMFs}\cup\{\bm I\}$
		\STATE $\bm r=\bm r-\bm I$
		\ENDWHILE
		%\STATE ${\rm IMFs=IMFs}\cup\{r\}$
	\end{algorithmic}
	\label{alg:FRIF}
\end{algorithm}

Notice that while the internal loop only consists of Hadamard multiplications among vectors, and its convergence properties can be analysed with the same tools used for the IF algorithm \cite{cicone2020iterative,cicone2020numerical}, in the outer loop we perform operations that may lead to a loss in accuracy of the method. We can thus adopt a spline interpolation to mitigate the accuracy loss, and even in this case, the computational cost of the outer loop is still $O(n\log n)$ operations due to the Fourier transforms.

As for the previous algorithms the matrix $K$, and thus the vector $\bm{\hat k_1}$, can be multiplied by a constant to upper bound its eigenvalues, and from Lemma \ref{res:IF_convergence} one can state an analogous convergence result.

\begin{corollary}
	Given a double-convoluted filter $k = \omega\star \omega$, then  the inner loop of the RIF Algorithm \ref{alg:RIFc} converges for any initial function $h(x)$.
	Moreover, if $M> 2$, then in the FRIF Algorithm \ref{alg:FRIF} the limit
	\[
	\lim_{m\to\infty}
	\bm{\hat h_m}
	=
	\lim_{m\to\infty}
	\bm{\hat k_1}^m \circ \bm{\hat h_1}
	\]
	converges for any vector $\bm{\hat h_1}$.
\end{corollary}

We have seen that the FRIF algorithm is provably convergent, and that its computational time is comparable with the FIF method. In the numerical examples, we will also show that empirically it produces sensible decompositions, but first let us address another property of the method.

\subsection{Anti-Aliasing Property}

In the discrete setting, the resampling of the signal $g(x)$ may in theory come with an undersampling of the highest frequencies, leading to aliasing effects. Here we show that in the FRIF algorithm, this is actually not a problem.

Suppose that the signal can be split into components $I_1(x), I_2(x), \dots$, where $I_1(x)$ has the highest instantaneous frequency among all the components.
In the FRIF algorithm we choose the resampling $G(x)$ where $ G^{-1}(z) = \int_{0}^x \ell(t)^{-1}\,{{\rm d}}t$ and $M=G^{-1}(1) = \int_{0}^1 \ell(t)^{-1}\,{{\rm d}}t$. The resampled signal $\wt h(x) = g(G(x))$ has thus domain $[0,M]$, but in the discrete setting we treat it as a signal over $[0,1]$, so we are actually working with
$$
h(x) = \wt h(Mx) = g(G(Mx)).
$$
The signal $h(x)$ presents now a new decomposition in components $J_1(x), J_2(x), \dots$, where $J_i(x) = I_i(G(Mx))$ and if $a_i(x)$ was the instantaneous frequency of $I_i(x)$, then the respective frequency of $J_i(x)$ is
$a_i(G(Mx)) G'(Mx) M$. Notice that the function $G(z)$ is chosen so that $J_1(x)$ is now approximately a stationary signal, so
\[
\beta \approx a_1(G(Mx)) G'(Mx) M = a_1(G(Mx))\ell(G(Mx))M,
\]
 for some constant $\beta$ and for every $x\in [0,1]$. As a consequence,
\[
a_1(G(Mx))G'(Mx)M  =
\int_{0}^1 \frac{a_1(G(Mx))G'(Mx)}{\ell(t)}
\,{{\rm d}}t	
\approx
\int_{0}^1 a_1(t)
\,{{\rm d}}t
\]
	that is surely less than $\|a_1(x)\|_\infty$.
Moreover, since $G(z)$ is increasing and $a_1(x) \ge a_i(x) \, \forall x,i$, then
\[
	a_1(G(Mx))G'(Mx)M\ge
	a_i(G(Mx))G'(Mx)M, \qquad \forall x
\]
	meaning that  $J_1(x)$  has still the biggest instantaneous frequency among the $J_i(x)$.	
This proves that the resampling does not create artificial high frequency components, so the
FRIF algorithm does not suffer from aliasing problems.

\subsection{Avoiding Interpolation}

As pointed out before, the interpolations may introduce a loss in accuracy on the output of Algorithm \ref{alg:FRIF}. One can though formulate a different, but equivalent, version of the continuous algorithm that does not require a resampling of the signal. Taking from the start of \eqref{eq:LALIF_as_first_order_RIF},
let $H(x) := G^{-1}(x)$, $x:=G(y)$ and $r:=G(y) - G(y-z)$, so that $z = y-G^{-1}(G(y) - r) = H(x) - H(x-r)$.
	\begin{align}\label{RIF2}
	\nonumber \int_{\f R} g(G(z))\,k(y-z) {{\rm d}}z &
	= \int_{\f R} g(G(y-z))\, k(z) {\rm d}z\\
	\nonumber &=  \int_{\f R} g(x - r)\, k(H(x) - H(x-r))
	H'(x-r) {\rm d}r
	\\&
	= \int_{\f R} g(r)\, k(H(x) - H(r))
	H'(r) {\rm d}r.
	\end{align}
	As a consequence, we can discretize the relation
	$$g_{m+1}(x)
	= g_m(x) - \int_{\f R} g_m(r)\, k(H(x) - H(r))
	H'(r) {\rm d}r$$
	by applying a quadrature rule on the points $x_i=i/(n-1)$, as
	\[
	g_{m+1}(x_i) \approx g_m(x_i) - \frac{1}{n-1}\sum_{j\in\mathbb Z}g_m(x_j) k(H(x_i)-H(x_j))H'(x_j)
	\]
	that coincides with multiplying the discretized signal $\bm g_m$ by the matrix $I_n - A_nD_n$, where
	\[
	A_n\hspace{-1pt}=\hspace{-1pt}\left[\frac{k(H(x_i)-H(x_j))}{n-1}\right]_{i,j=0}^{n-1}\hspace{-2pt}=\hspace{-1pt}\left[
	\frac{
		k\left(H\left(\frac{i}{n-1}\right)-H\left(\frac j{n-1}\right)\right)
	}{n-1}
	\right]_{i,j=0}^{n-1},
	\qquad D_n = \diag\left(H'(x_i)_{i=0}^{n-1}\right).
	\]
Notice that $D_n$ is positive definite, since from \eqref{eq:resampling_from_length}, $H'(x) = \ell(x)^{-1}>0$ and $A_n$ is symmetric since the filter $k$ is an even function.
If we call $B_N$ the matrix $A_N$ in the case $H(x)\equiv x$, then by Corollary 3 of \cite{cicone2020numerical}, $B_N$ is positive semidefinite. Since for a big enough $N$, the matrix $A_n$ is approximated up to an arbitrary small error by a $n\times n$ principal submatrix of the matrix $B_N$, then we can conclude that $A_n$ is also positive semidefinite. As a consequence,
\[
I_n - A_nD_n \sim
I_n -  D_n^{1/2}A_n D_n^{1/2}
\]
and all its eigenvalues are real and less than 1. Eventually, as in the precedent algorithms, the matrix $A_nD_n$ can be multiplied by a constant so that its eigenvalues are upper bounded for example by 1, so that $\Lambda(I_n-A_nD_n)\cu (-1,1]$ and the method becomes provably convergent.

The resulting algorithm is thus equivalent in its continuous version to Algorithm \ref{alg:RIFc}, and in its discrete version it avoids the need to interpolate the signal two times per IMF. Moreover, its internal loop has been proved to be convergent and it presents the same flexibility properties as ALIF.

At the same time, though, the matrix $A_nD_n$ is not cyclic, so we lose the fast implementation that was possible in Algorithm \ref{alg:FRIF}.
%Moreover, usually we don't know the length function $\ell(x)$ and thus the inverse sampling $H(x)$ in their analytic form, but only through their value on a grid, so the exact computation of $A_n$ and $D_n$ may still introduce errors.
For this reason, we do not test this version of the RIF algorithms in the following numerical experiments.

\section{Numerical Experiments}\label{sec:NumericalExamples}

In this section we show and compare the performances of all the reviewed techniques. In order to study the signals and their decompositions in time-frequency, we will rely on the so called IMFogram, a recently developed algorithm \cite{Barbe2021time}, which allows to represent the frequency content of all IMFs. The IMFogram proves to be a robust, fast and reliable way to obtain the time-frequency representation of a signal, and it has been shown to converge, in the limit, to the well know spectrogram based on the FFT \cite{cicone2021IMFogram}.

The following tests have been conducted using MATLAB$^{\tiny{\textregistered}}$ R2021a installed on a 64--bit Windows 10 Pro computer equipped with a 11th Gen Intel$^{\tiny{\textcopyright}}$
 Core$^{\tiny{\textregistered}}$  i7-1165G7 at 2.80GHz processor and 8GB RAM. All tested examples and algorithms are freely available at \footnote{\url{www.cicone.com}}.

\subsection{Example 1}

\begin{figure}
\includegraphics[width=0.33\linewidth]{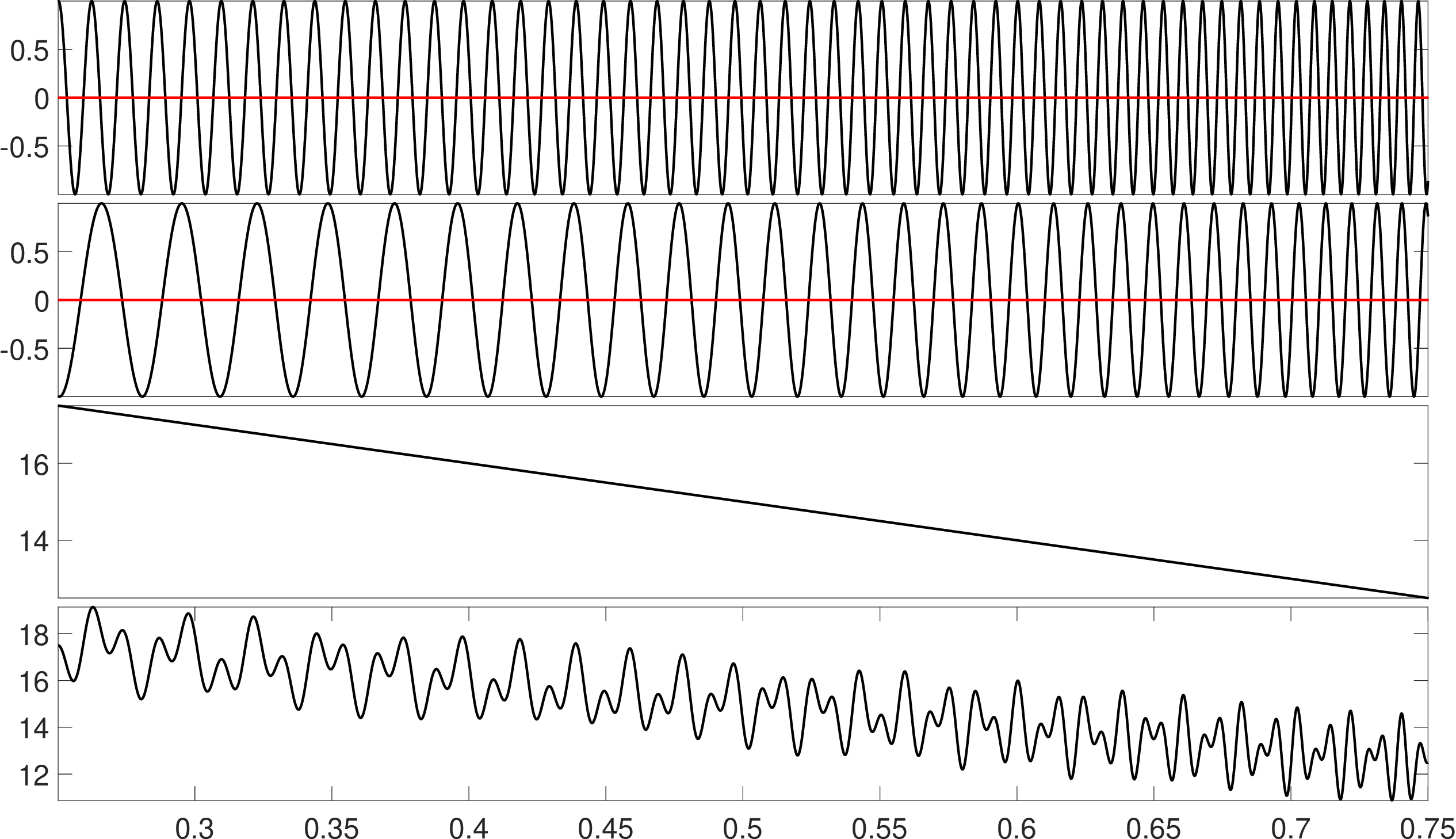}~\includegraphics[width=0.33\linewidth]{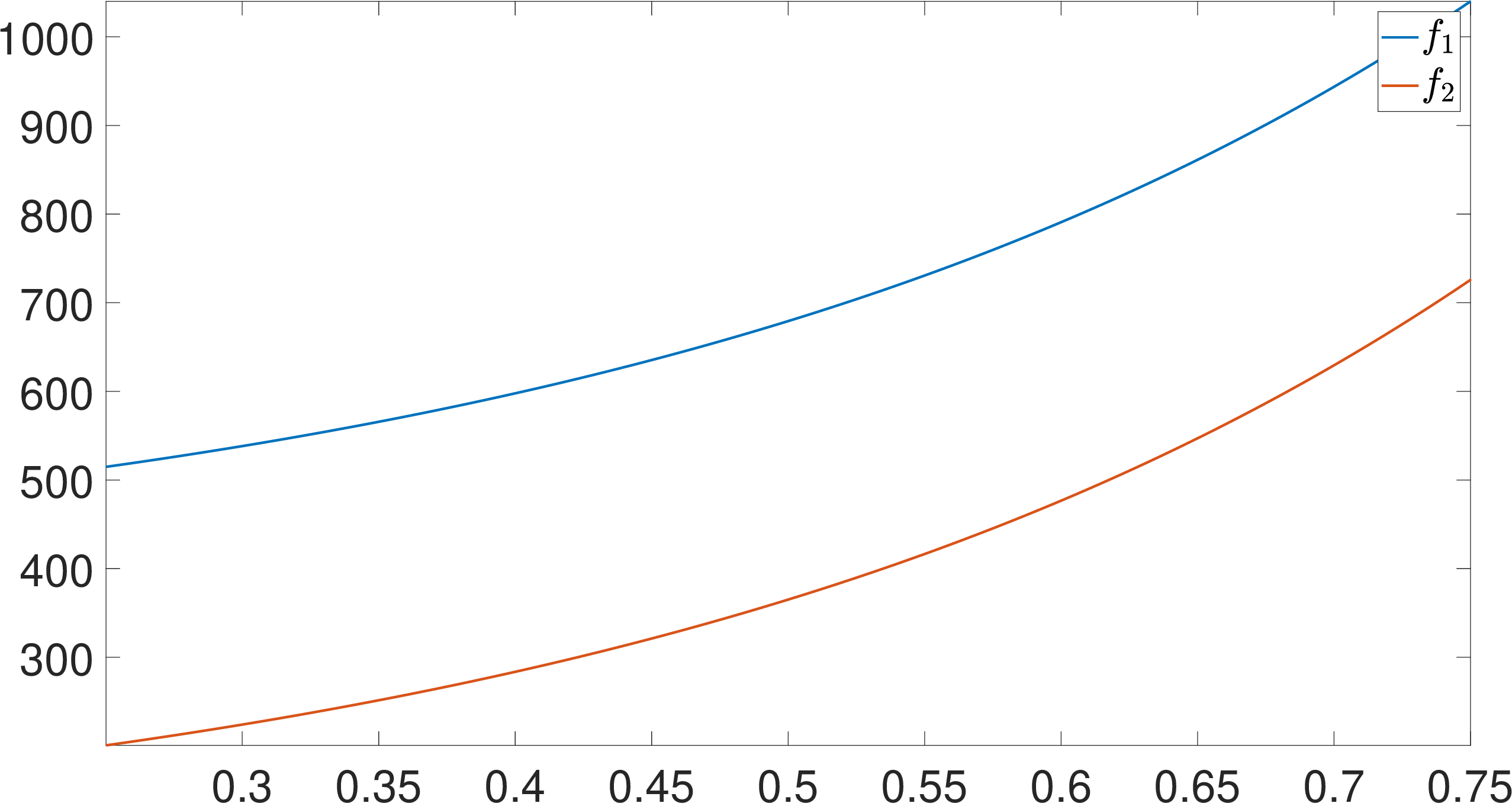}~\includegraphics[width=0.33\linewidth]{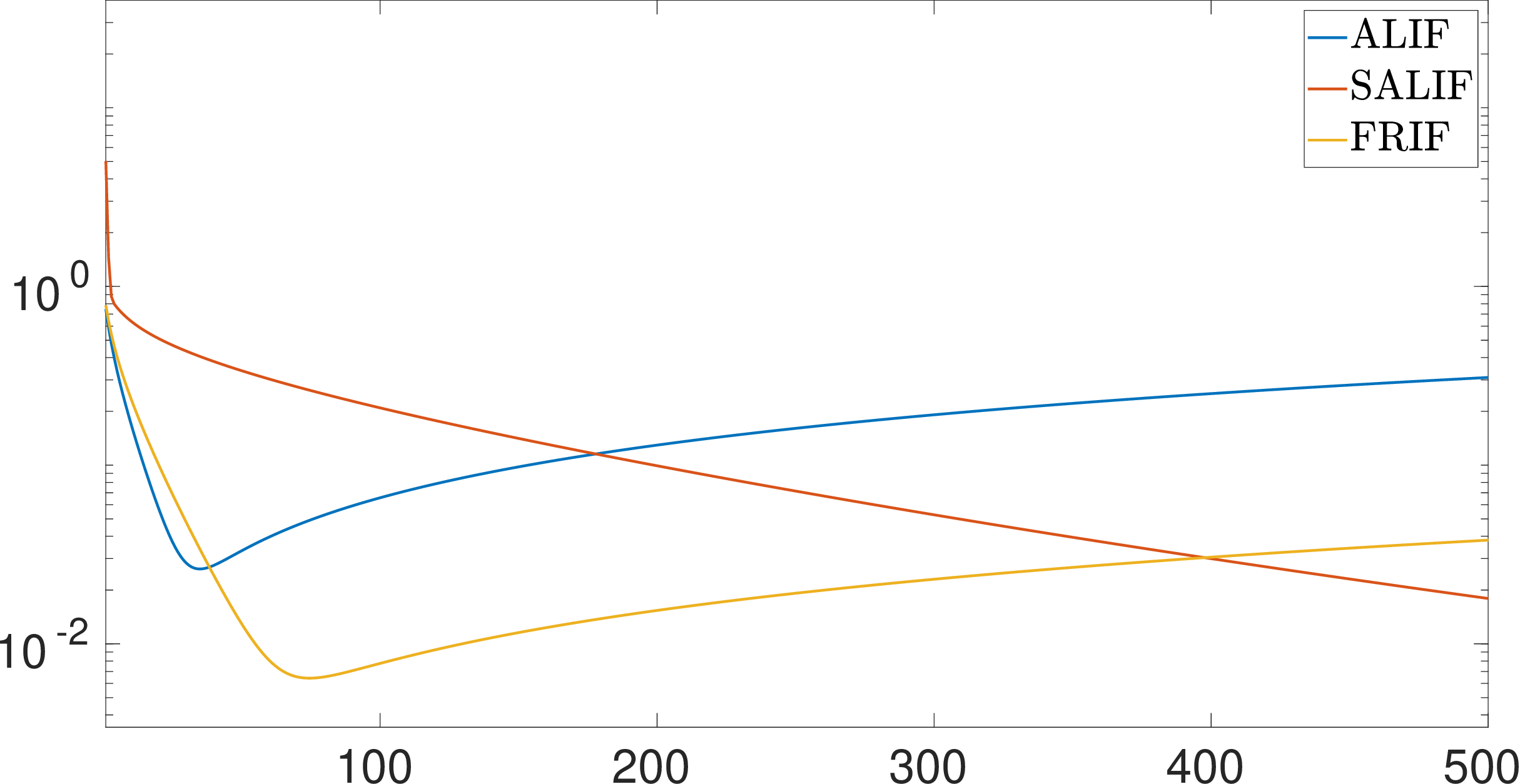}
\caption{Example 1. Left panel: the components $f_1$ and $f_2$, respectively first and second row,the trend, third row, and the signal $f$, bottom row. Central panel: exponential instantaneous frequencies of $f_1$ and $f_2$. Right panel: relative error in norm 2 between the ground truth and $\IMF_1$ produced by ALIF, SALIF, and FRIF algorithms.}
\label{fig:Ex1_sig}
\end{figure}

We consider the artificial signal $f$, plotted in the left panel, bottom row, of Figure \ref{fig:Ex1_sig}, which contains two nonstationary components with exponentially changing frequencies $f_1$ and $f_2$, plus a trend $f_3$. In particular

\begin{eqnarray}
  f_1(x) &=& \cos(20 e^{t\pi}+120\pi t) \\
  f_2(x) &=& \cos(20 e^{t\pi}+20\pi t) \\
  f_3(x) &=& -10x+20 \\
\end{eqnarray}
where $x$ vary in $[0,\, 1]$ and is sampled over $10^4$ points.

The $f_1$ and $f_2$ components and $f$ signal are plotted in the left panel of Figure \ref{fig:Ex1_sig}, whereas $f_1$ and $f_2$ frequencies are shown in the central panel.

%\begin{figure}
%\includegraphics[width=0.33\linewidth]{Ex1_IMFs_ALIF_minus_GT.eps}~\includegraphics[width=0.33\linewidth]{Ex1_IMFs_SALIF_minus_GT.eps}~\includegraphics[width=0.33\linewidth]{Ex1_IMFs_FRIF_minus_GT.eps}
%\caption{Example 1. Differences between the ground truth and the decompositions derived using ALIF (left), SALIF (central), FRIF (right).}
%\label{fig:Ex1_IMFs_minus_GT}
%\end{figure}

in Table \ref{tab:Ex1} we report the computational time required by ALIF, SALIF and FRIF with a fixed stopping criterion.
%when the parameter delta in the stopping criterion is set to $10^{-6}$.
In the same table we summarize the performance of the three techniques in terms of inner loop iterations required to produce the two IMFs and the relative error measured as ratio between the norm 2 of the difference between the computed IMF and the corresponding ground truth, and the norm 2 of the ground truth itself.

\begin{table}[h!]
\centering
\begin{tabular}{|c||c|c|c|}
  \hline
  Example 1 & ALIF & SALIF & FRIF  \\
     \hline
       \hline
time(s)  & $16.3293$ & $26.9395$ & $\bm{0.9107}$ \\
\hline
$\textrm{err}_1$ & $0.040260$ & $0.117824$ & $\bm{0.006535}$ \\
\hline
$\textrm{err}_2$ & $1.051461$ & $0.117842$ & $\bm{0.006543}$ \\
\hline
$\textrm{err}_3$ & $0.049352$ & $0.000084$ & $\bm{0.000017}$ \\
\hline
num of iter $\textrm{IMF}_1$  & $61$ & $175$ & $80$ \\
\hline
num of iter $\textrm{IMF}_2$  & $500$ & $155$ & $ 4$ \\
  \hline
\end{tabular}
\caption{performance of various techniques when applied on Example 1, measured as relative errors in norm 2 and number of iterations.}
\label{tab:Ex1}
\end{table}

From Table \ref{tab:Ex1} results %, as well as from Figure \ref{fig:Ex1_IMFs_minus_GT} plots,
it is clear that FRIF proves to converge quickly to a really accurate solution. In fact, it takes less than a second to produce a decomposition which has a relative error which is order of magnitudes smaller than the ones produced using ALIF and SALIF methods. Furthermore ALIF and SALIF decompositions require more than 16 and 26 seconds, respectively, to converge. This is confirmed by the results shown in the right panel of Figure \ref{fig:Ex1_sig}, where we compare the norm 2 relative error of the $\IMF_1$ obtained using ALIF, SALIF, and FRIF algorithms for subsequent steps in the inner loops when we remove the stopping condition.
%, minus the corresponding ground truth, and the norm 2 of the ground truth itself.
ALIF initially tends toward the right solution. At 35 steps the relative error reach the minimum value of $0.0262$, and then, after that, the instabilities of the method show up and drive the solution far away from the right one. SALIF, instead, is clearly convergent, in fact the solution is moving steadily to the exact one. However SALIF converge rate is small, as proven by the relative error which is slowly decaying. In fact, after 500 inner loop steps, the relative error is still around $0.0179$.
Finally, FRIF quickly converge to a really good approximation of the right solution, at 73 steps the error is minimal with a relative error of $0.0064$. After this step, the relative error restarts growing due to the chosen stopping criterion. It is important to remember, in fact, that, in general, the ground truth is not known. This is the reason why the stopping criterion adopted in these techniques does not rely on the ground truth knowledge. Hence, as a consequence, FRIF, ALIF and SALIF, do not necessarily stop when the actual best approximation of the ground truth is achieved. For example, one can see that the ALIF algorithm doesn't stop in the computation of the second IMF of the signal.
Studying what it could be an ideal stopping criterion and how to tune it properly is outside the scope of this work.
%In this example, for comparison reasons, we simply set the stopping criterion so that all techniques perform at least 500 iterations.

\subsection{Example 2}

In this second example, we start from the artificial signal $h$ which contains two nonstationary components, $h_1$ and $h_2$, and a trend $h_3$,
%then we add three different realizations of gaussian noise.

\begin{eqnarray}
  h_1(x) &=& \cos(20\cos(4\pi t)-160\pi t) \\
  h_2(x) &=& \cos(20\cos(4\pi t)-280\pi t) \\
  h_3(x) &=& \cos(2\pi t) \\
\end{eqnarray}
where $x$ vary in $[0,\, 1]$ and is sampled over 8000 points.

\begin{figure}
\includegraphics[width=0.49\linewidth]{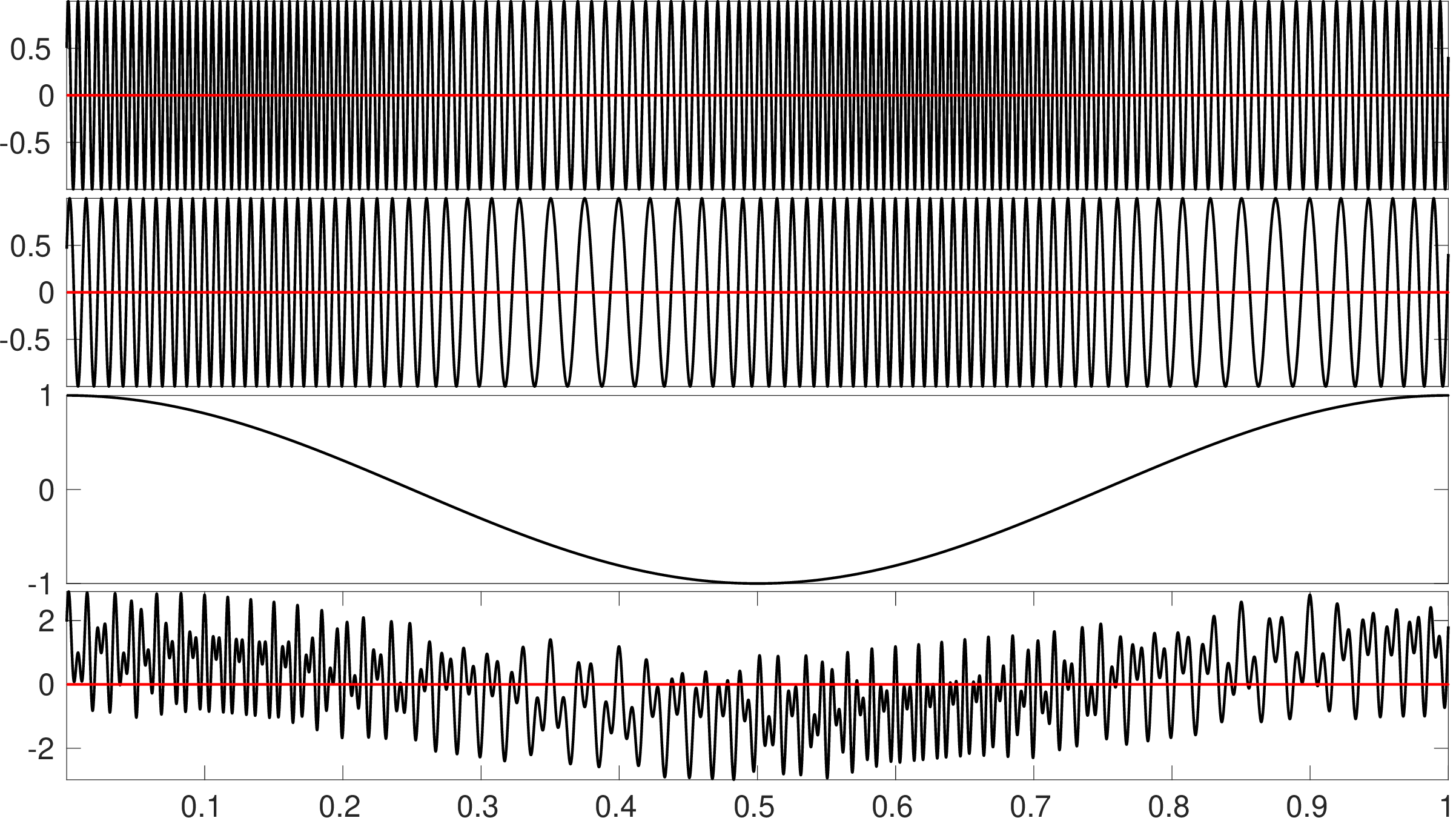}~\includegraphics[width=0.49\linewidth]{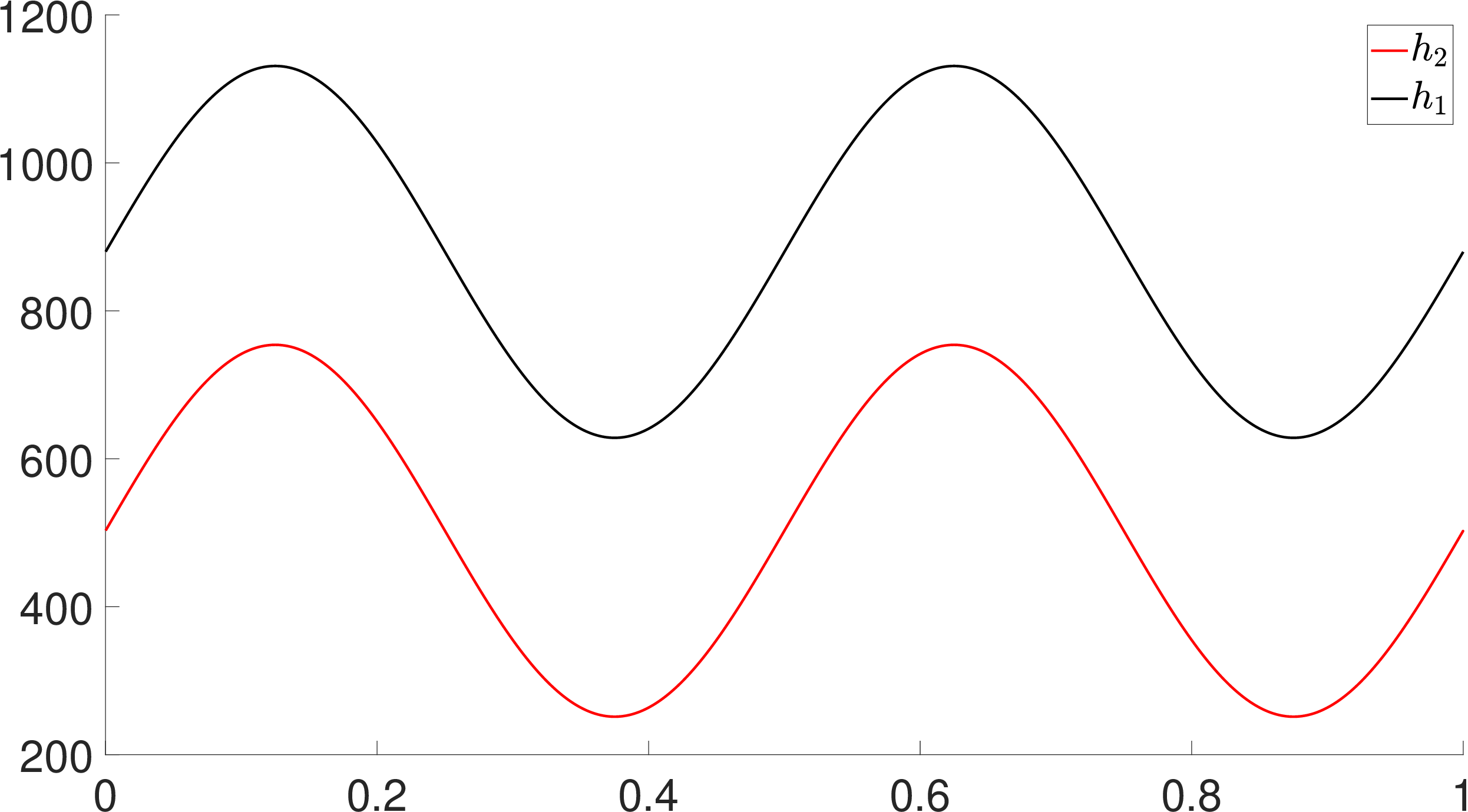}
\caption{Example 2. Left panel: the components $h_1$ and $h_2$, respectively first and second row, and the signal $h$, bottom row. Right panel: exponential instantaneous frequencies of $h_1$ and $h_2$.}\label{fig:Ex2_sig}
\end{figure}

The $h_1$, $h_2$, the trend component, and $h$ signal are plotted in the left column of Figure \ref{fig:Ex2_sig}, whereas $h_1$ and $h_2$ frequencies are shown in the right panel.

\begin{figure}
\includegraphics[width=0.33\linewidth]{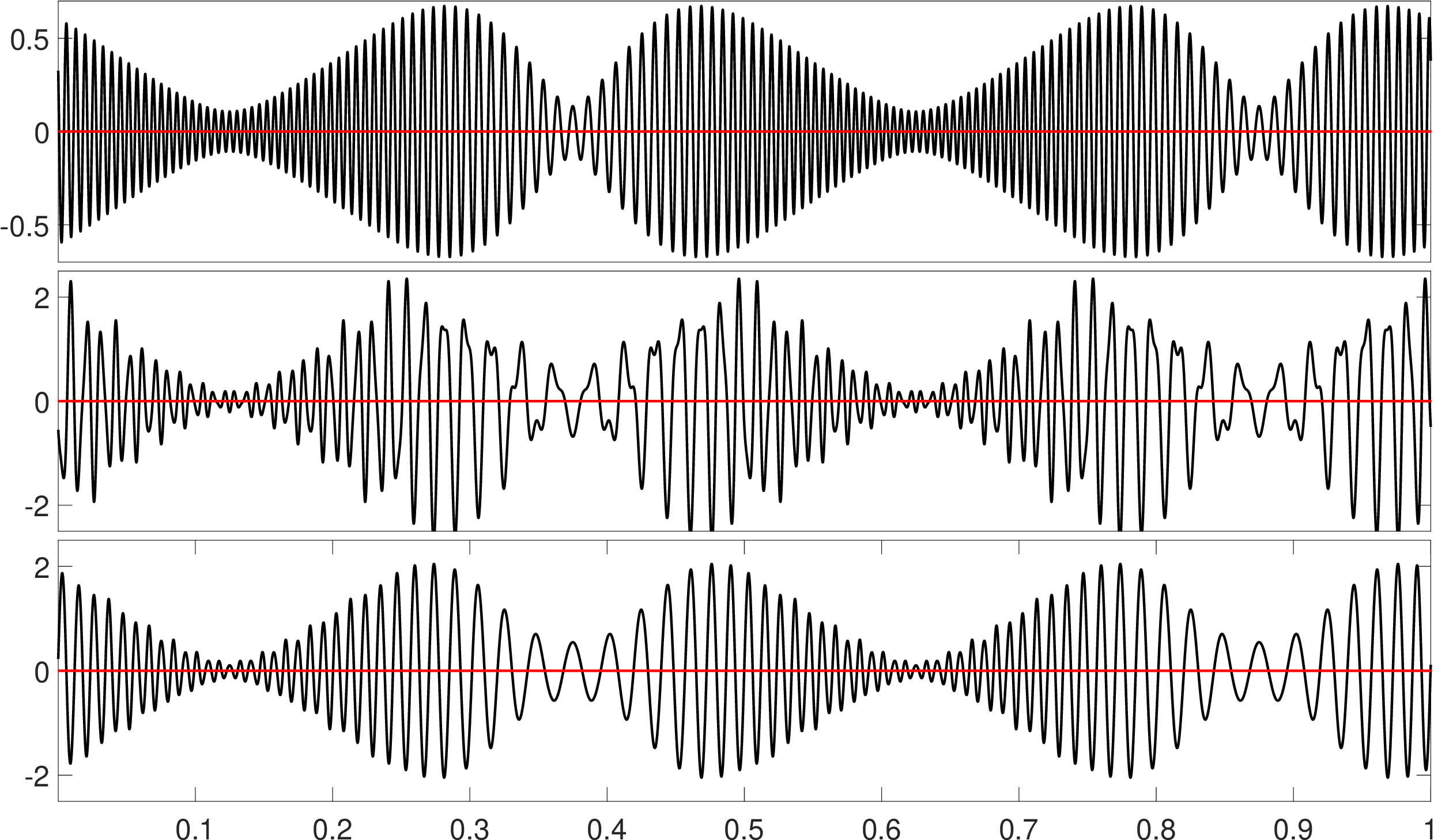}~\includegraphics[width=0.33\linewidth]{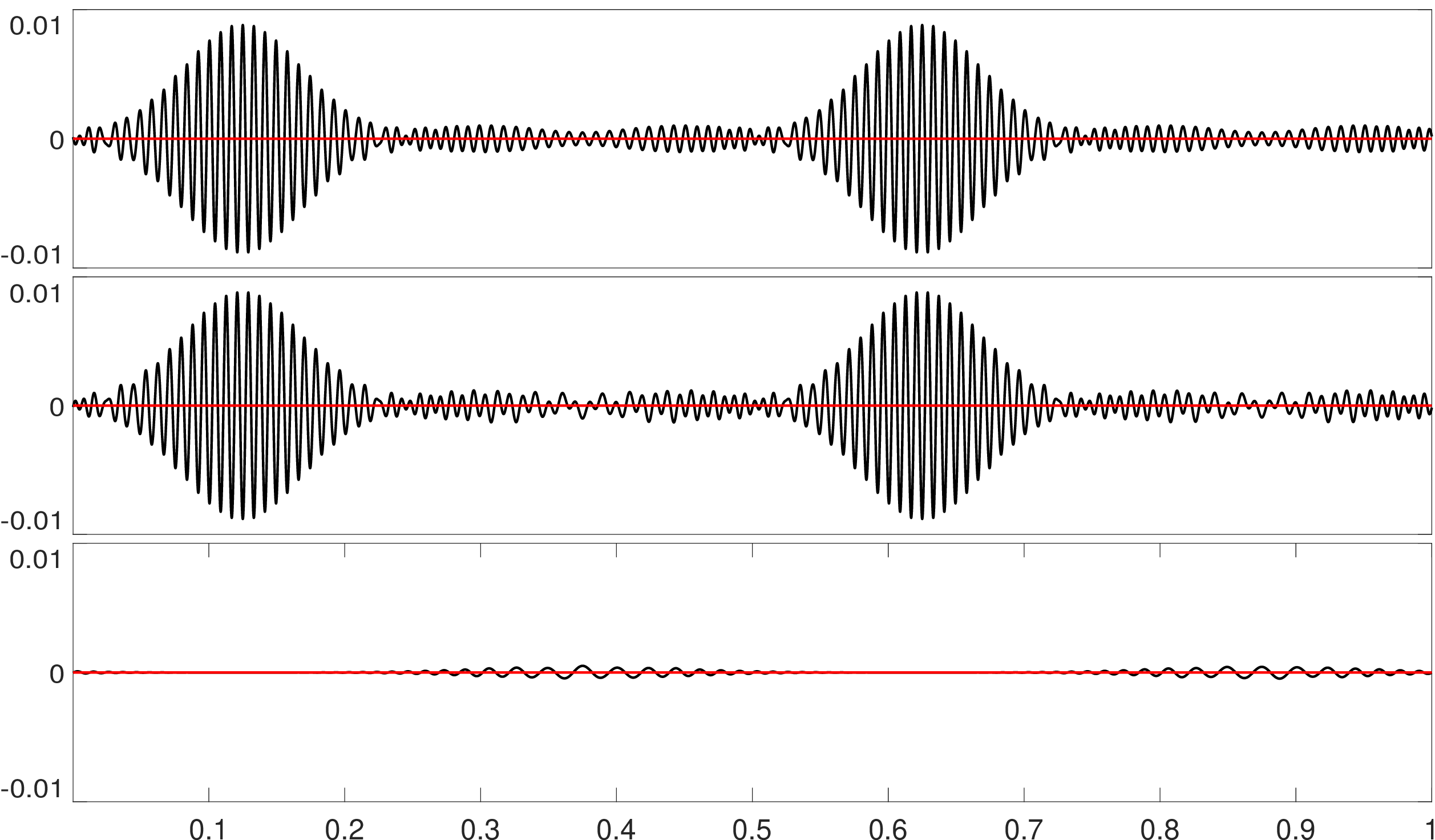}~\includegraphics[width=0.33\linewidth]{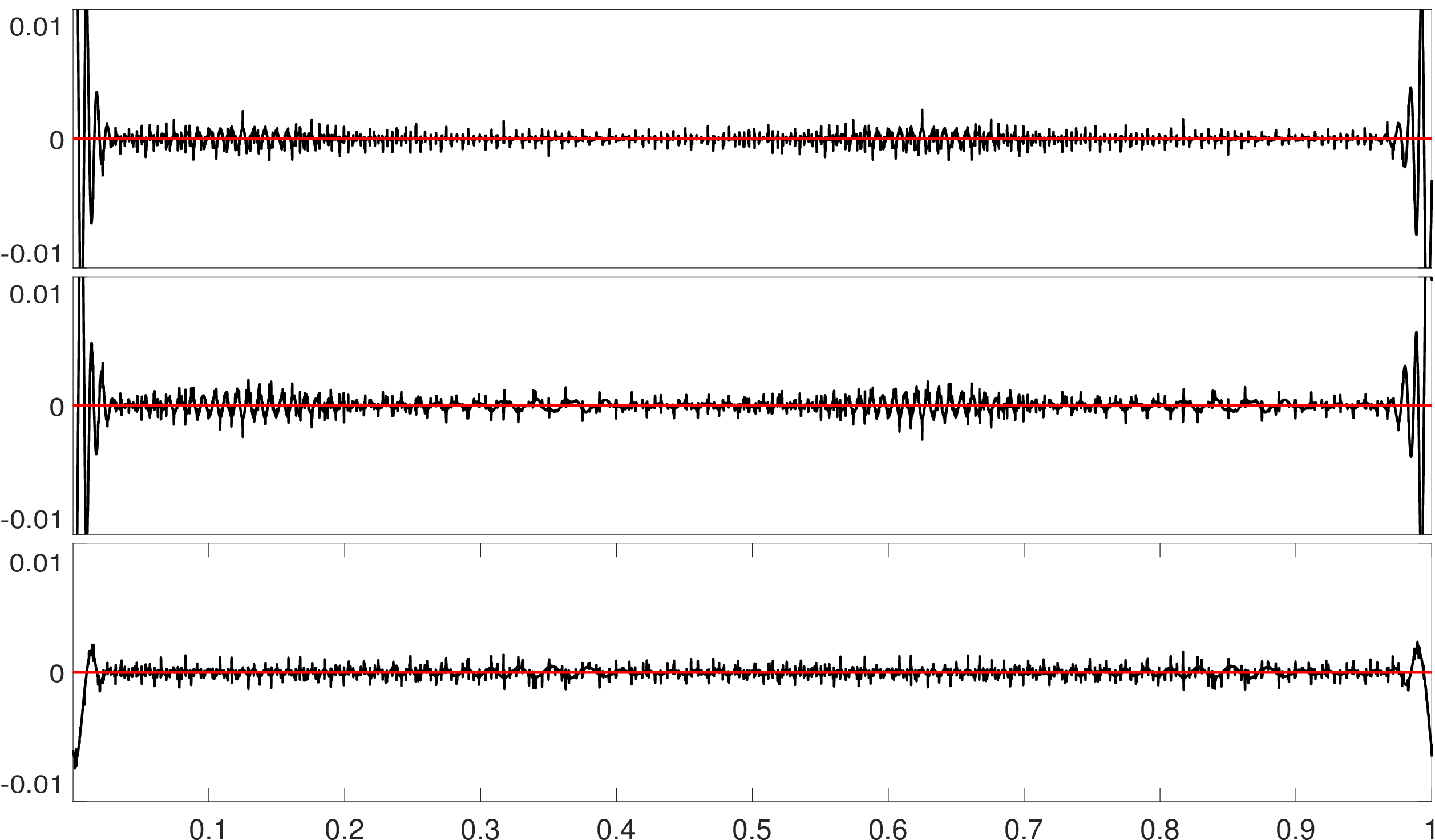}
\caption{Example 2. Difference between the ground truth and the derived decomposition via ALIF (left), SALIF (central), FRIF (right).}\label{fig:Ex2_IMFs}
\end{figure}

in Table \ref{tab:Ex2} we report the performance of ALIF, SALIF and FRIF techniques. In Figure \ref{fig:Ex2_IMFs} we show the differences between the IMFs produced by the different methods and the known ground truth. It is evident both from the table and the figure that the proposed FRIF method outperform the other approaches both from the computational and the accuracy point of view.

\begin{table}[h!]
\centering
\begin{tabular}{|c||c|c|c|}

  \hline
  Example 2 & ALIF & SALIF & FRIF  \\
     \hline
       \hline
time(s)  & $16.0005$ & $20.2120$ & $\bm{1.0958}$ \\
\hline
$\textrm{err}_1$ & $0.457672$ & $0.003584$ & $\bm{0.003426}$ \\
\hline
$\textrm{err}_2$ & $1.374017$ & $0.003591$ & $\bm{0.003292}$ \\
\hline
$\textrm{err}_3$ & $1.304946$ & $\bm{0.000229}$ & $0.000908$ \\
\hline
num of iter $\textrm{IMF}_1$  & $500$ & $468$ & $81$ \\
\hline
num of iter $\textrm{IMF}_2$  & $500$ & $ 8$ & $11$ \\
  \hline
\end{tabular}
\caption{Example 2 performance of ALIF, SALIF and FRIF, measured as relative errors in norm 2 and iteration number.}
\label{tab:Ex2}
\end{table}

\subsection{Example 3}

\begin{figure}
\includegraphics[width=0.45\linewidth]{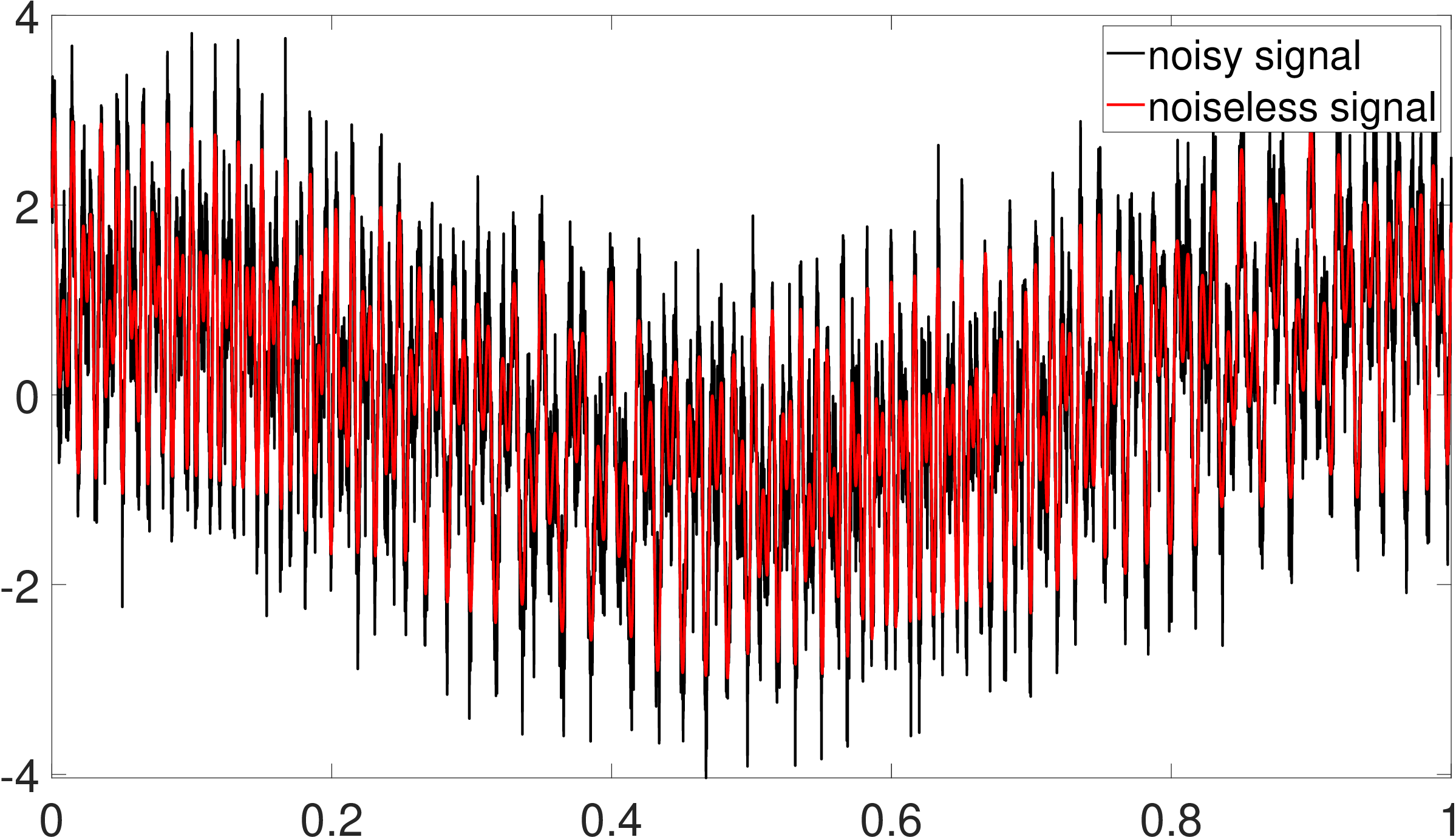}~\includegraphics[width=0.45\linewidth]{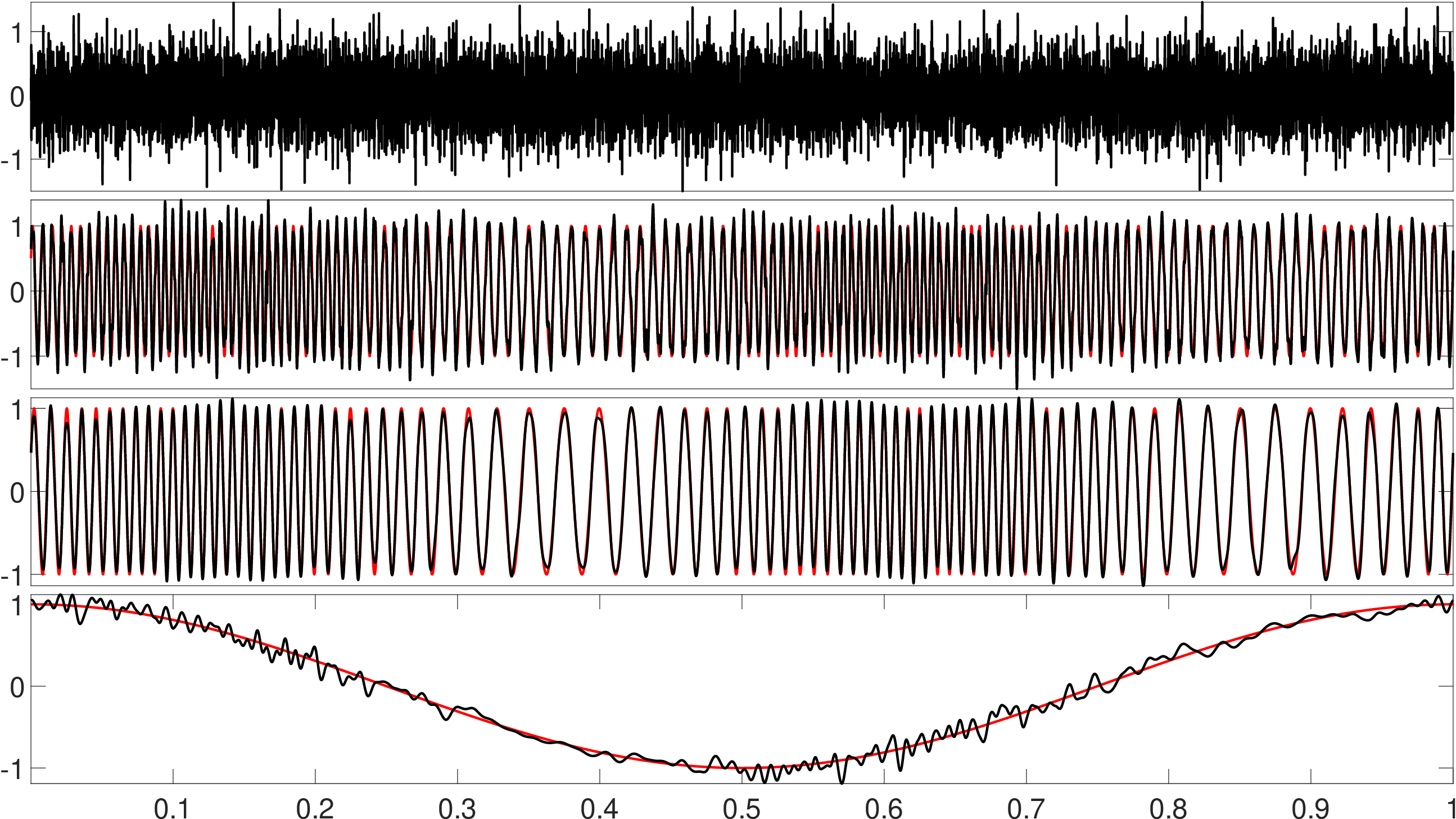}
\caption{Example 3. Left panel, the noisy signal compared with the noiseless signal $h$ defined in Example 2. The SNR is around 8.6 dB. Right panel, the IMF decomposition derived by FRIF.}\label{fig:Ex3_1}
\end{figure}

In this example we show the robustness of the proposed FRIF approach to noise. To do so, we consider the signal $h$ studied in Example 2 and we perturb it by additive Gaussian noise. In Figure \ref{fig:Ex3_1} we plot on the left panel the perturbed signal when the signal to noise ratio (SNR) is of 8.6 dB. On the right panel we report the decomposition produced by FRIF. It is evident that the method can separate properly the random perturbation in the first row, from the deterministic components in the following three rows.

\begin{figure}
\includegraphics[width=0.45\linewidth]{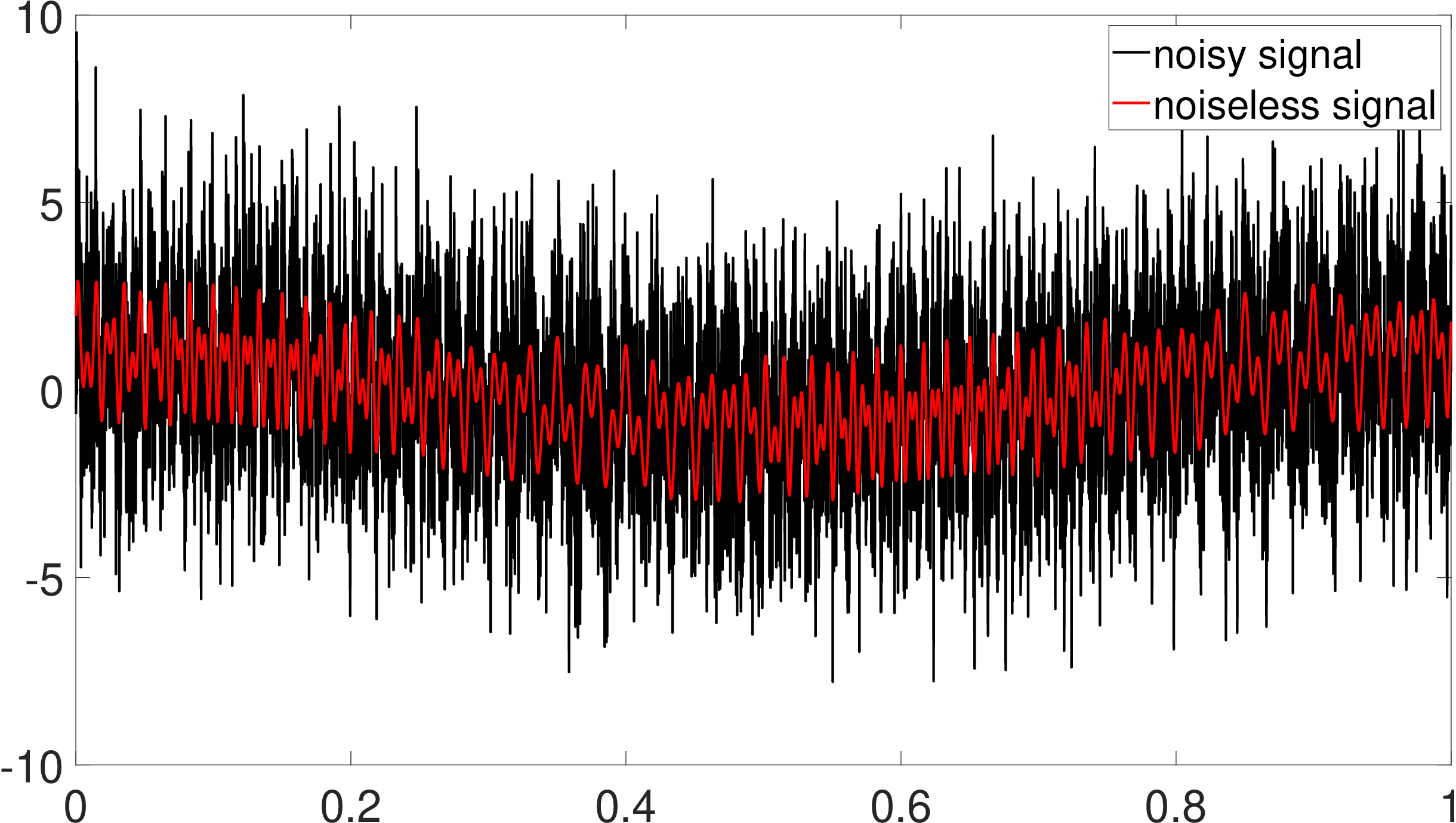}~\includegraphics[width=0.45\linewidth]{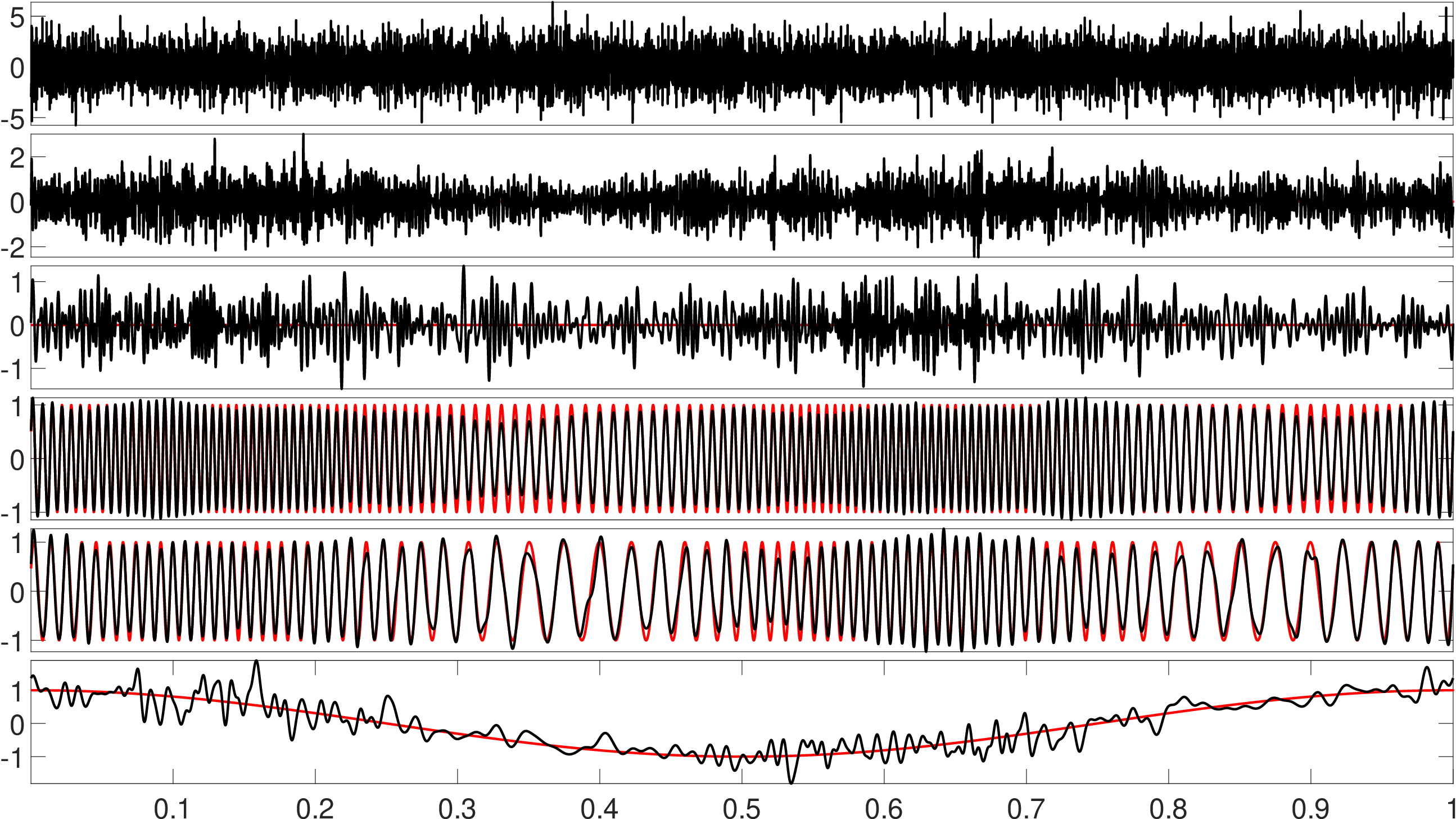}
\caption{Example 3. Left panel, the noisy signal with SNR around 1.3 dB compared with the noiseless signal $h$ of Example 2. Right panel, the corresponding FRIF decomposition compared with the ground truth.}\label{fig:Ex3_2}
\end{figure}

 This result is confirmed even if we increase the SNR to 1.3 dB, left panel of Figure \ref{fig:Ex3_2}. It is evident from this figure that this level of noise is quite high. Nevertheless FRIF method proves to be able still to separate the deterministic signal from the additive Gaussian contribution, as shown in the left panel of Figure \ref{fig:Ex3_2}.

\subsection{Example 4}

\begin{figure}
\includegraphics[width=0.31\linewidth]{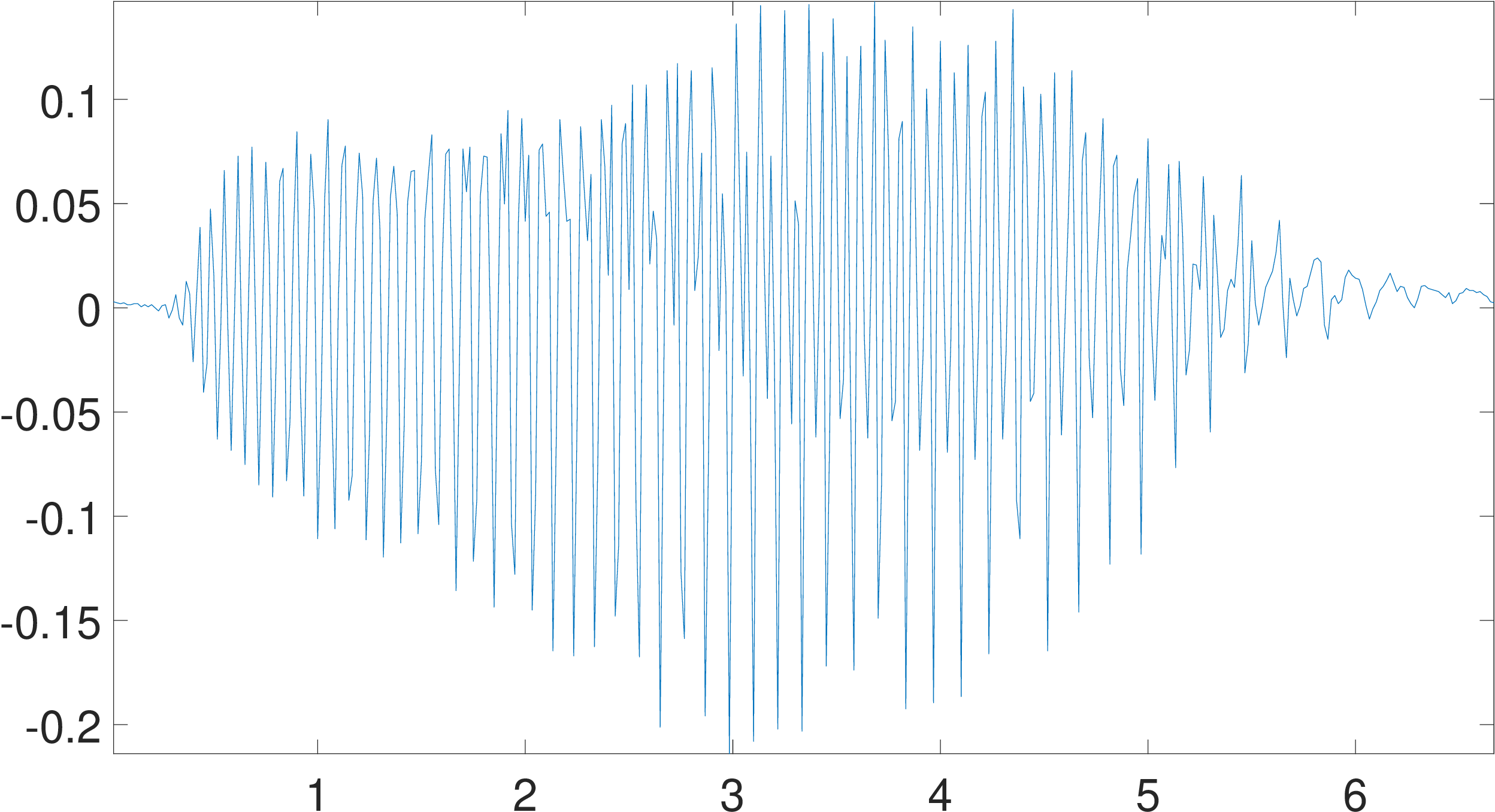}~\includegraphics[width=0.36\linewidth]{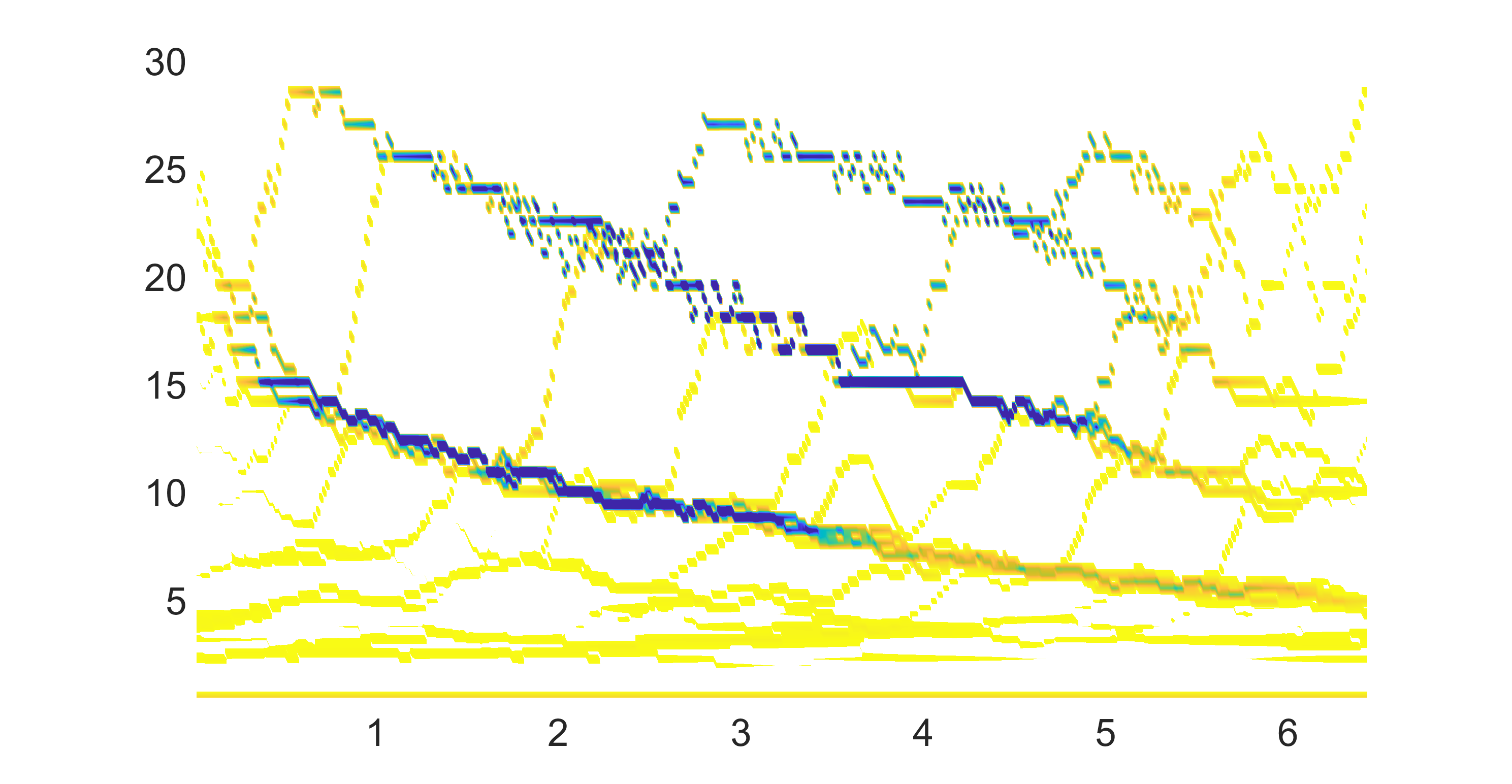}~\includegraphics[width=0.31\linewidth]{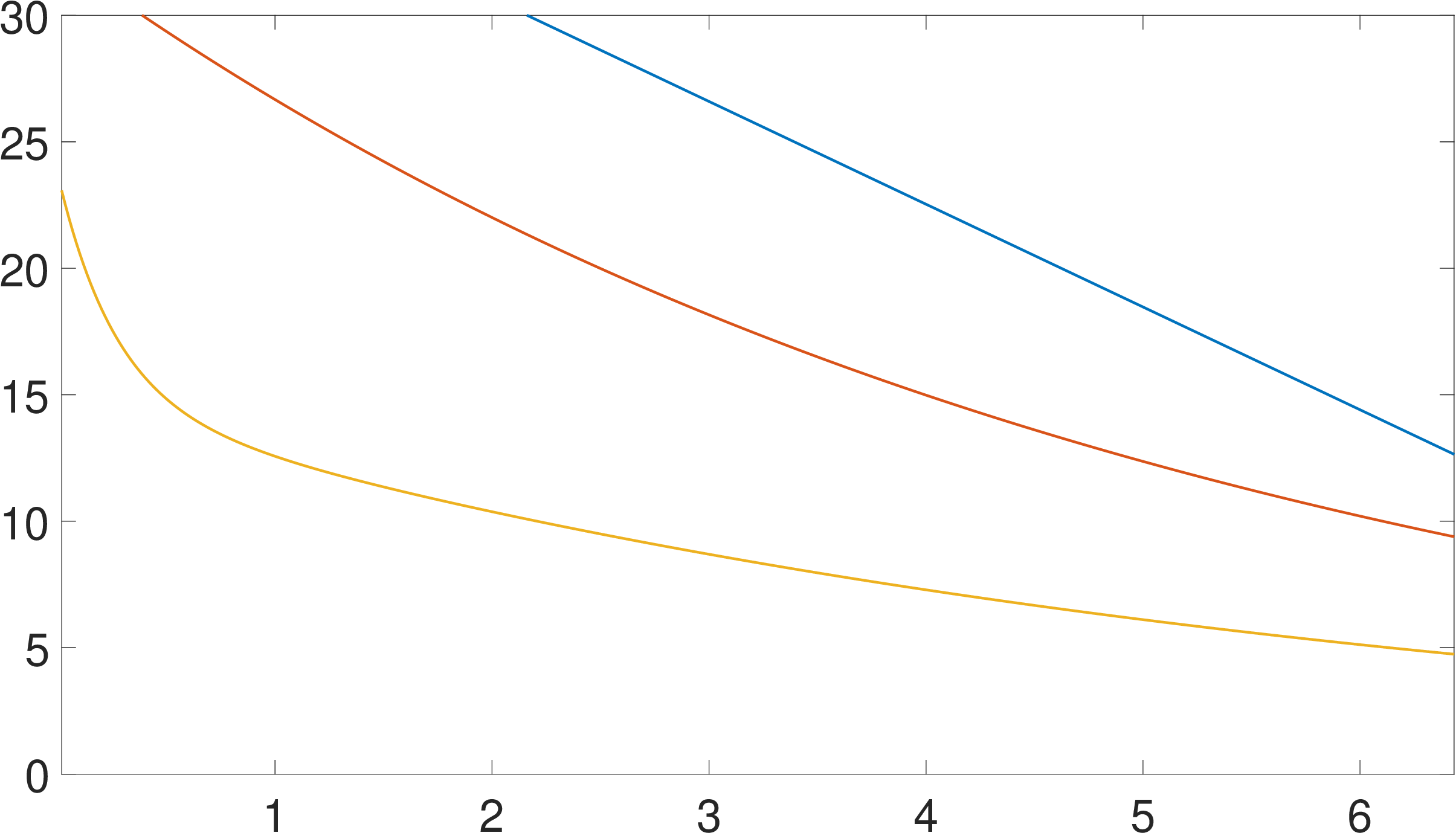}
\caption{Example 4. Left panel, sound produced by a bat. Central panel, the corresponding IMFogram time-frequency plot. Right panel, instantaneous frequency curves inferred from the IMFogram plot.}\label{fig:Ex4_sig}
\end{figure}

We conclude the numerical section with an example based on a real life signal. We consider the recording of the sound emitted by a bat, shown in the left panel of Figure \ref{fig:Ex4_sig}. In the central panel, we show the associated time-frequency plot obtained using the IMFogram \cite{Barbe2021time}. From this plot we observe that this signal appears to contain three main simple oscillatory components which present rapid changes in frequencies. Those are classical examples of the so called chirps. By using a curve extraction method, it is possible to derive from the IMFogram the instantaneous frequency curves plotted in the right panel of Figure \ref{fig:Ex4_sig}. As briefly mentioned earlier, the identification of these instantaneous frequency curves is of fundamental importance for the proper functioning of FRIF, but it is also a research topic per se. In this work, we assume that they can be computed accurately and we postpone the analysis of how to compute them in a robust and accurate way to future works.

\begin{figure}
\includegraphics[width=0.25\linewidth]{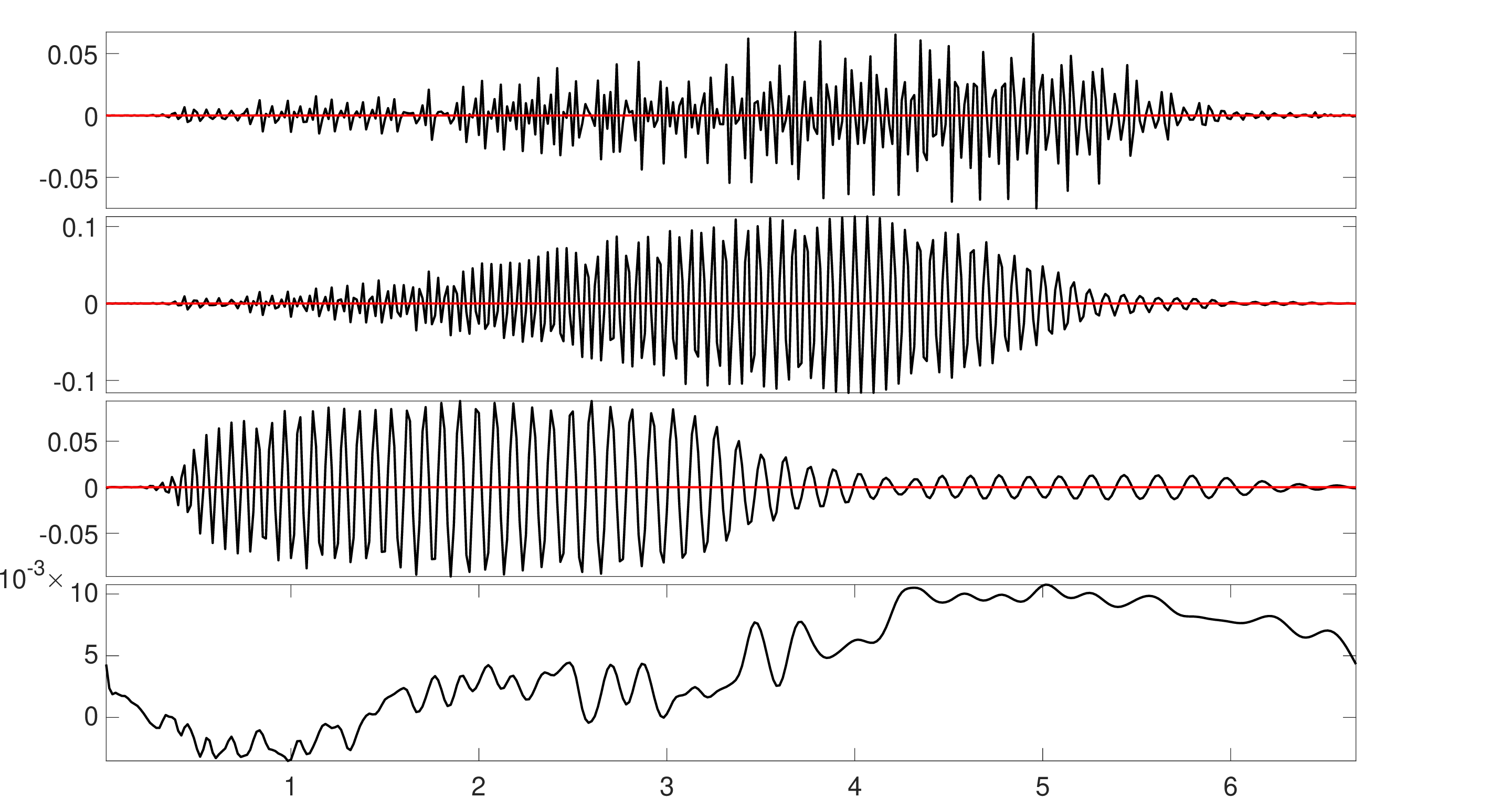}~\includegraphics[width=0.25\linewidth]{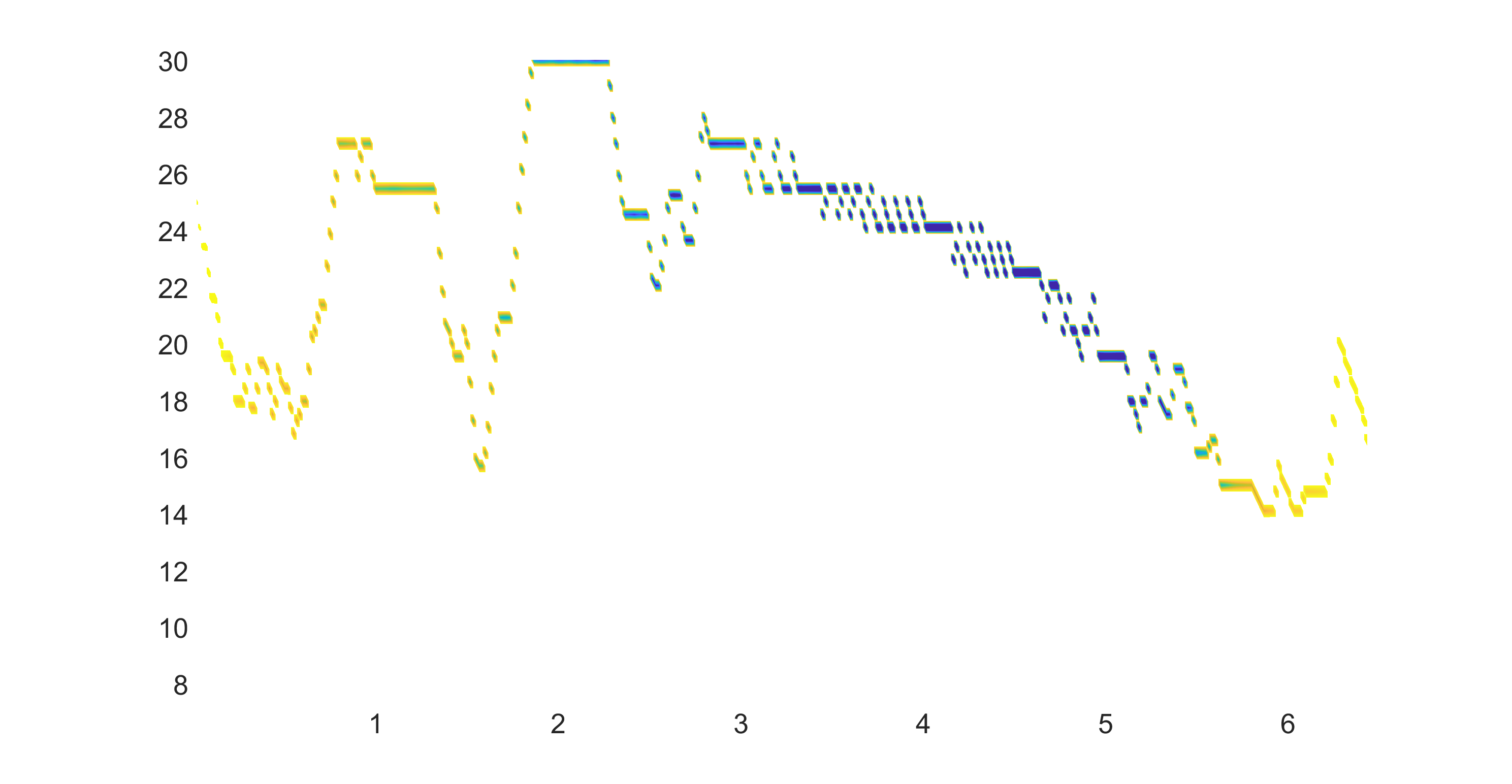}~\includegraphics[width=0.25\linewidth]{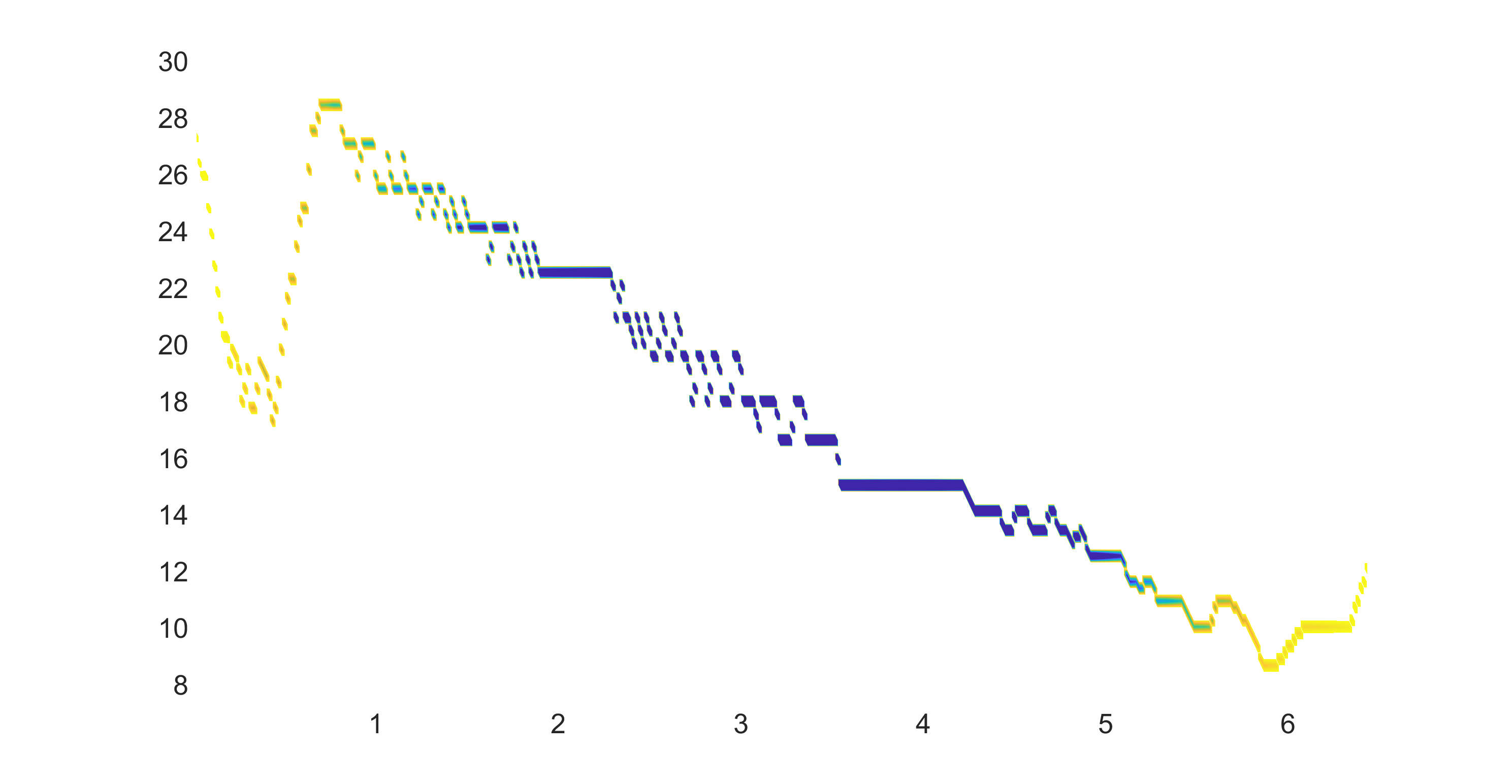}~\includegraphics[width=0.25\linewidth]{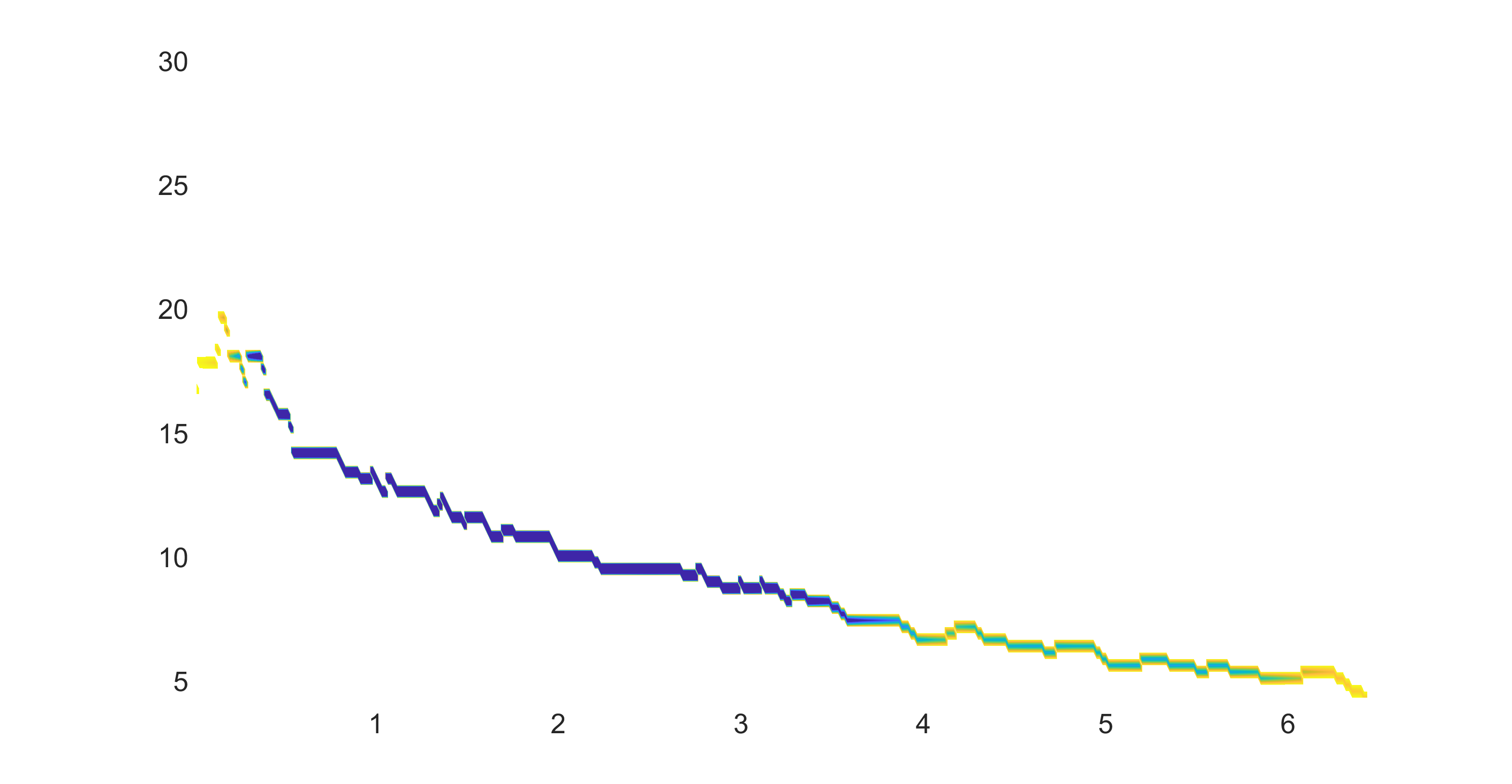}
\caption{Example 4. Left most panel, IMF decomposition produced by FRIF. From central left to right most panel, the IMFogram time-frequency plots associated with the first, second and third row in the FRIF decomposition, respectively.}\label{fig:Ex4_IMFs}
\end{figure}

By leveraging on the extracted curves, we run FRIF algorithm and derive the decomposition shown in the left most panel of Figure \ref{fig:Ex4_IMFs}. The first three IMFs produced correspond to the three main chirps observed in the IMFogram, which is depicted in the central panel of Figure \ref{fig:Ex4_sig}. This is confirmed by running IMFogram separately on the first three IMFs produced by FRIF. The results are shown in the rightmost 3 panels of Figure \ref{fig:Ex4_IMFs}. From these plots it becomes clear that the algorithm is able to separate in a clean way the three chirps contained in the signal.

\section{Conclusions}\label{sec:Outlook}

Following the success of the Empirical Mode Decomposition (EMD) method for the decomposition of non-stationary signals, and given that its mathematical understanding is still very limited, in recent years the Iterative Filtering (IF) first, and then the Adaptive Local Iterative Filtering (ALIF) have been proposed. They inherit the same structure of EMD, but rely on convolution for the computation of the signal moving average. On the one hand, the mathematical understanding of IF is now pretty advanced, this include its acceleration in what is called Fast Iterative Filtering (FIF) and its complete convergence analysis. On the other hand, IF proved to be limited in separating, in a physically meaningful way, components which exhibit quick changes in their frequencies, like chirps or whistles. For this reason ALIF was proposed as a generalization of IF which overcome the limitations that are present in IF. However, even though some advances have been obtained in recent years, the theoretical understanding of ALIF is far from being complete. In particular, it is not yet clear under which assumptions it is possible to guarantee a priori its convergence.

For this reason, in this work we introduced the Resampled Iterative Filtering (RIF), and, in the discrete setting, the Stable Adaptive Local Iterative Filtering (SALIF) and the Fast Resampled Iterative Filtering (FRIF), that are capable of decomposing non-stationary signals into simple oscillatory components, even in presence of fast changes in their instantaneous frequencies, like in chirps. We have analyzed them from a theoretical stand point, showing, among other things, that it is possible to guarantee a priori their convergence. Furthermore, we have tested them using several artificial and real life examples.

More is yet to be said about the argument. In particular, all these methods are dependent on the computation of a length function $\ell(x)$ which is, de facto, the reciprocal of the instantaneous frequency curve associated with each component contained in the signal. This function is required to guide the aforementioned methods, including ALIF itself, in the extraction of physically meaningful IMFs. The identification of instantaneous frequency curves associated with each component contained in a given signal is a research topic per se, and it is out of the scope of the present work. This is why we plan to study this problem in future works.

Another open problem regards the selection of an optimal stopping criterion and its tuning to be used in this kind of methods. The stopping criterion implemented can influence consistently the performance of these techniques. We plan to work in this direction in the future.

Finally, we plan to work on the extension of the proposed techniques to handle multidimensional and multivariate signals.

\section*{Acknowledgements}

The authors are members of the Italian ``Gruppo Nazionale di Calcolo Scientifico'' (GNCS) of the Istituto Nazionale di Alta Matematica ``Francesco Severi'' (INdAM). 

AC thanks the Italian Ministry of the University and Research for the financial support under the PRIN PNRR 2022 grant number E53D23018040001 ERC field PE1 project P2022XME5P titled ``Circular Economy from the Mathematics for Signal Processing prospective'', and the Ministry of Foreign Affairs and the International Cooperation for the financial support under the ``Grande Rilevanza'' Italy – China Science and Technology Cooperation Joint Project titled ``sCHans – Solar loading infrared thermography and deep learning teCHniques for the noninvAsive iNSpection of precious artifacts''.

GB is supported by the Alfred Kordelinin s\"a\"ati\"o grant no. 210122 and partly by 
an Academy of Finland grant (Suomen Akatemian p\"a\"at\"os 331240). GB is also supported by the European Union (ERC consolidator grant, eLinoR, no 101085607).

All authors contributed to the work in equal parts. All authors read and approved the final manuscript.\\

This preprint has not undergone peer review or any post-submission improvements or corrections. The Version of Record of this article is published in Numerische Mathematik and is available online at\\ https://doi.org/10.1007/s00211-024-01394-y. 

\bibliographystyle{abbrv}

\begin{thebibliography}{10}

\bibitem{An2017}
X.~An.
\newblock Local rub-impact fault diagnosis of a rotor system based on adaptive
  local iterative filtering.
\newblock {\em Transactions of the Institute of Measurement and Control},
  39(5):748--753, 2017.

\bibitem{an2016application}
X.~An, C.~Li, and F.~Zhang.
\newblock Application of adaptive local iterative filtering and approximate
  entropy to vibration signal denoising of hydropower unit.
\newblock {\em Journal of Vibroengineering}, 18(7):4299--4311, 2016.

\bibitem{an2016wind}
X.~An and L.~Pan.
\newblock Wind turbine bearing fault diagnosis based on adaptive local
  iterative filtering and approximate entropy.
\newblock {\em Proceedings of the Institution of Mechanical Engineers, Part C:
  Journal of Mechanical Engineering Science}, 231(17):3228--3237, 2017.

\bibitem{an2017vibration}
X.~An, W.~Yang, and X.~An.
\newblock Vibration signal analysis of a hydropower unit based on adaptive
  local iterative filtering.
\newblock {\em Proceedings of the Institution of Mechanical Engineers, Part C:
  Journal of Mechanical Engineering Science}, 231(7):1339--1353, 2017.

\bibitem{an2016demodulation}
X.~An, H.~Zeng, and C.~Li.
\newblock Demodulation analysis based on adaptive local iterative filtering for
  bearing fault diagnosis.
\newblock {\em Measurement}, 94:554--560, 2016.

\bibitem{barbarino2021conjectures}
G.~Barbarino and A.~Cicone.
\newblock Conjectures on spectral properties of alif algorithm, 2021.
\newblock arXiv:2009.00582.

\bibitem{Barbe2021time}
P.~Barbe, A.~Cicone, W.~Suet~Li, and H.~Zhou.
\newblock Time-frequency representation of nonstationary signals: the imfogram.
\newblock {\em Pure and Applied Functional Analysis}, 2021.

\bibitem{cicone2020iterative}
A.~Cicone.
\newblock Iterative filtering as a direct method for the decomposition of
  nonstationary signals.
\newblock {\em Numerical Algorithms}, pages 1--17, 2020.

\bibitem{cicone2020study}
A.~Cicone and P.~Dell'Acqua.
\newblock Study of boundary conditions in the iterative filtering method for
  the decomposition of nonstationary signals.
\newblock {\em Journal of Computational and Applied Mathematics}, 373:112248,
  2020.

\bibitem{cicone2019spectral}
A.~Cicone, C.~Garoni, and S.~Serra-Capizzano.
\newblock Spectral and convergence analysis of the discrete alif method.
\newblock {\em Linear Algebra and its Applications}, 580:62--95, 2019.

\bibitem{cicone2021IMFogram}
A.~Cicone, W.~S. Li, and H.~Zhou.
\newblock New theoretical insights in the decomposition and time-frequency
  representation of nonstationary signals: the imfogram algorithm.
\newblock {\em preprint}, 2021.

\bibitem{cicone2016adaptive}
A.~Cicone, J.~Liu, and H.~Zhou.
\newblock Adaptive local iterative filtering for signal decomposition and
  instantaneous frequency analysis.
\newblock {\em Applied and Computational Harmonic Analysis}, 41(2):384--411,
  2016.

\bibitem{cicone2020convergence}
A.~Cicone and H.-T. Wu.
\newblock Convergence analysis of adaptive locally iterative filtering and sift
  method.
\newblock {\em submitted}, 2021.

\bibitem{cicone2020numerical}
A.~Cicone and H.~Zhou.
\newblock Numerical analysis for iterative filtering with new efficient
  implementations based on fft.
\newblock {\em Numerische Mathematik}, 147(1):1--28, 2021.

\bibitem{huang2009convergence}
C.~Huang, L.~Yang, and Y.~Wang.
\newblock Convergence of a convolution-filtering-based algorithm for empirical
  mode decomposition.
\newblock {\em Advances in Adaptive Data Analysis}, 1(04):561--571, 2009.

\bibitem{huang2014introduction}
N.~E. Huang.
\newblock Introduction to the hilbert--huang transform and its related
  mathematical problems.
\newblock {\em Hilbert--Huang transform and its applications}, pages 1--26,
  2014.

\bibitem{huang1998empirical}
N.~E. Huang, Z.~Shen, S.~R. Long, M.~C. Wu, H.~H. Shih, Q.~Zheng, N.-C. Yen,
  C.~C. Tung, and H.~H. Liu.
\newblock The empirical mode decomposition and the hilbert spectrum for
  nonlinear and non-stationary time series analysis.
\newblock {\em Proceedings of the Royal Society of London. Series A:
  mathematical, physical and engineering sciences}, 454(1971):903--995, 1998.

\bibitem{kim2016multiscale}
S.~J. Kim and H.~Zhou.
\newblock A multiscale computation for highly oscillatory dynamical systems
  using empirical mode decomposition (emd)--type methods.
\newblock {\em Multiscale Modeling \& Simulation}, 14(1):534--557, 2016.

\bibitem{li2018entropy}
Y.~Li, X.~Wang, Z.~Liu, X.~Liang, and S.~Si.
\newblock The entropy algorithm and its variants in the fault diagnosis of
  rotating machinery: A review.
\newblock {\em IEEE Access}, 6:66723--66741, 2018.

\bibitem{lin2009iterative}
L.~Lin, Y.~Wang, and H.~Zhou.
\newblock Iterative filtering as an alternative algorithm for empirical mode
  decomposition.
\newblock {\em Advances in Adaptive Data Analysis}, 1(04):543--560, 2009.

\bibitem{mitiche2018classification}
I.~Mitiche, G.~Morison, A.~Nesbitt, M.~Hughes-Narborough, B.~G. Stewart, and
  P.~Boreham.
\newblock Classification of partial discharge signals by combining adaptive
  local iterative filtering and entropy features.
\newblock {\em Sensors}, 18(2):406, 2018.

\bibitem{cicone2017Geophysics}
M.~Piersanti, M.~Materassi, A.~Cicone, L.~Spogli, H.~Zhou, and R.~G. Ezquer.
\newblock Adaptive local iterative filtering: A promising technique for the
  analysis of nonstationary signals.
\newblock {\em Journal of Geophysical Research: Space Physics},
  123(1):1031--1046, 2018.

\bibitem{sharma2017automatic}
R.~Sharma, R.~B. Pachori, and A.~Upadhyay.
\newblock Automatic sleep stages classification based on iterative filtering of
  electroencephalogram signals.
\newblock {\em Neural Computing and Applications}, 28(10):2959--2978, 2017.

\bibitem{stallone2020new}
A.~Stallone, A.~Cicone, and M.~Materassi.
\newblock New insights and best practices for the successful use of empirical
  mode decomposition, iterative filtering and derived algorithms.
\newblock {\em Scientific Reports}, 10:15161, 2020.

\bibitem{torres2011complete}
M.~E. Torres, M.~A. Colominas, G.~Schlotthauer, and P.~Flandrin.
\newblock A complete ensemble empirical mode decomposition with adaptive noise.
\newblock In {\em 2011 IEEE international conference on acoustics, speech and
  signal processing (ICASSP)}, pages 4144--4147. IEEE, 2011.

\bibitem{ur2011filter}
N.~Ur~Rehman and D.~P. Mandic.
\newblock Filter bank property of multivariate empirical mode decomposition.
\newblock {\em IEEE transactions on signal processing}, 59(5):2421--2426, 2011.

\bibitem{WU2020current}
H.-T. Wu.
\newblock Current state of nonlinear-type time-frequency analysis and
  applications to high-frequency biomedical signals.
\newblock {\em Current Opinion in Systems Biology}, 23:8--21, 2020.

\bibitem{wu2009ensemble}
Z.~Wu and N.~E. Huang.
\newblock Ensemble empirical mode decomposition: a noise-assisted data analysis
  method.
\newblock {\em Advances in adaptive data analysis}, 1(01):1--41, 2009.

\bibitem{yang2017oscillation}
D.~Yang, B.~Wang, G.~Cai, and J.~Wen.
\newblock Oscillation mode analysis for power grids using adaptive local
  iterative filter decomposition.
\newblock {\em International Journal of Electrical Power \& Energy Systems},
  92:25--33, 2017.

\bibitem{yeh2010complementary}
J.-R. Yeh, J.-S. Shieh, and N.~E. Huang.
\newblock Complementary ensemble empirical mode decomposition: A novel noise
  enhanced data analysis method.
\newblock {\em Advances in adaptive data analysis}, 2(02):135--156, 2010.

\bibitem{zheng2014partly}
J.~Zheng, J.~Cheng, and Y.~Yang.
\newblock Partly ensemble empirical mode decomposition: An improved
  noise-assisted method for eliminating mode mixing.
\newblock {\em Signal Processing}, 96:362--374, 2014.

\end{thebibliography}

\end{document}